\documentclass[11 pt]{smfart}
\usepackage[latin1]{inputenc}
\usepackage[T1]{fontenc}
\usepackage[english,french]{babel}
\usepackage{ae, aecompl}
\usepackage{amsmath,amsfonts,amssymb}
\usepackage{amsthm}
\usepackage{mathrsfs,array,graphicx}
\usepackage{stmaryrd}
\usepackage{wasysym}
\usepackage[all,ps]{xy}
\usepackage{calc}
\usepackage{enumerate}
\usepackage{ulem}
\usepackage{smfthm}
\usepackage{smfenum}

\usepackage{fancyhdr}
\addto\captionsfrench{\def\proofname{Preuve}}

\newcommand{\uni}{\varpi_F}
\newcommand{\ugu}{\stackrel{\mathrm{d\acute{e}f}}{=}}

\newcommand{\Q}{\mathbb{Q}_p}

\newcommand{\Z}{\mathbb{Z}_p}

\newcommand{\OF}{\mathcal{O}_F}
\newcommand{\Oe}{\mathcal{O}_E}

\newcommand{\C}{C}

\newcommand{\Cr}{C^{r}}

\newcommand{\Ind}{\mathrm{Ind}}
\newcommand{\into}{\hookrightarrow}

\newcommand{\GL}{\mathrm{GL}_2}

\newcommand{\val}{\mathrm{val}}

\newcommand{\Gal}{\mathrm{Gal}}

\newcommand{\Sym}{\mathrm{Sym}}

\numberwithin{equation}{section}

\DeclareMathOperator{\Hom}{Hom}
\input xy
\xyoption{all}

\usepackage{amscd}

\author[M. De Ieso]{Marco De Ieso}
\address{Bâtiment 430, Université Paris-Sud, 91405, Orsay Cedex, France}

\email{Marco.DeIeso@math.u-psud.fr}

\title{SUR CERTAINS COMPLÉTÉS UNITAIRES UNIVERSELS EXPLICITES POUR $\GL(F)$}

\begin{document}
\frontmatter
\date{}



\pagestyle{fancy}






\begin{abstract}
Dans cet article, nous donnons une description explicite du complété unitaire universel de certaines représentations localement $\Q$-analytiques de $\GL(F)$, où $F$ est une extension finie de $\Q$ (ce qui généralise des résultats de Berger-Breuil pour $F = \Q$). Pour cela, nous utilisons 
certains espaces de Banach de fonctions de classe $C^r$ sur $\OF$ (pour $r$ dans $\mathbb{R}_{\geq 0}$) introduits dans \cite{marco}.
\end{abstract}

\begin{altabstract}
In this paper we give an explicit description of the universal unitary completion of certain locally $\Q$-analytic representations of $\GL(F)$, where $F$ is a finite extension of $\Q$ (this generalizes some results of Berger-Breuil for $F=\Q$). To this aim, we make use of certain Banach spaces of $C^r$ functions on $\OF$ (for $r\in \mathbb{R}_{\geq 0}$) introduced in \cite{marco}.
\end{altabstract}


\maketitle

\tableofcontents

\section{Introduction} 

Soit $p$ un nombre premier. La dernière décennie a vu l'émergence et la preuve d'une correspondance locale $p$-adique entre certaines représentations continues de dimension $2$ de $\Gal(\overline{\mathbb{Q}}_p/\Q)$ et certaines représentations de $\GL(\Q)$. Cette correspondance, qui a pris le nom de correspondance de Langlands $p$-adique pour $\GL(\Q)$, a été initiée par  Breuil (\cite{breuilab}, \cite{breuila}) et a été établie par  Colmez \cite{colmez2} et  Pa\v{s}k\={u}nas \cite{pask} à la suite de travaux de Colmez \cite{colmez3} et Berger-Breuil \cite{bb}. 

Si $F$ est une extension finie de $\Q$, $F\neq \Q$, la question d'associer des représentations $p$-adiques de $G\ugu \GL(F)$ aux représentations $p$-adiques de dimension $2$ de $\Gal(\overline{\mathbb{Q}}_p/F)$ dans l'esprit d'une correspondance locale à la Langlands n'est  pas encore comprise et les résultats obtenus pour l'instant sont très partiels. Cependant Breuil \cite{bre}, en utilisant principalement les travaux de Schraen \cite{scr} et de Frommer \cite{fr} sur la filtration de Jordan-Hölder des induites paraboliques localement $\Q$-analytiques, définit une représentation localement $\Q$-analytique $\Pi(V)$ de $G$ pour la plupart des représentations cristallines $V$ de dimension $2$ de $\Gal(\overline{\mathbb{Q}}_p/F)$ et à poids de Hodge-Tate distincts, et en commence l'étude. En général, la représentation $\Pi(V)$ ne permet pas de reconstruire la représentation galoisienne de départ, toutefois on s'attend à ce qu'elle intervienne  comme sous-objet de la bonne représentation, ce qui fait que les complétés unitaires universels  de ses constituants fondamentaux sont des objets pertinents.

L'objet du présent article est celui de donner une description explicite du complété unitaire universel de certaines induites paraboliques localement $\Q$-analytiques (en particulier celles qui interviennent dans la construction de la représentation $\Pi(V)$). La motivation du fait qu'une  telle description est possible est suggérée  par \cite[Theorème 4.3.1]{bb}, où les auteurs décrivent le complété unitaire universel d'une induite parabolique localement algébrique de $\GL(\Q)$ en utilisant l'espace des fonctions de classe $C^r$ sur $\Z$, $r$ étant un nombre rationnel positif qui dépend de l'induite considerée.  

Pour cela l'auteur a introduit et exploré dans \cite{marco} une nouvelle notion de fonction de classe $C^r$ sur $\OF$, où $r$ désigne un nombre réel positif et $\OF$ l'anneau d'entiers de $F$, qui s'appuie principalement sur les travaux d'Amice, Amice-Velù, Vishik, Van der Put et Colmez (\cite{ami}, \cite{amivel}, \cite{vis}, \cite{vander}, \cite{colmez}) et qui repose sur l'idée  cruciale suivante: une fonction $f$  de $\OF$ dans $E$ est de classe $C^r$ si $f(x+y)$ a un développement limité à l'ordre $[r]$ (où $[r]$ désigne la partie entière de $r$) en tout $x$, et si le reste est $o(|y|^r)$ uniformément (en $x$) sur tout compact. Dans \cite{marco} on a aussi montré que cette notion ne coïncide pas avec une autre définition naturelle de fonction de classe $C^r$ sur $\OF$ obtenue en voyant $\OF$ comme $\Z^{[F:\Q]}$ (Remarque \ref{alternativa}).

Voir si les complétés unitaires universels que nous avons construits sont non nuls est, en général,  une question délicate  et complétement résolue seulement dans le cas $F=\Q$ \cite[Corollaire 5.3.1]{bb} en utilisant la théorie des $(\varphi,\Gamma)$-modules de Fontaine \cite{font}. Mentionnons par ailleurs que le \cite[Theorème 4.3.1]{bb} est un ingrédient important pour établir ce résultat. En déhors de $\Q$ nous ne connaissons pas en général la réponse. Toutefois, on déduit la non nullité dans quelques cas à partir des résultats de Vigneras \cite{vig} (voir aussi \cite{ks} pour une preuve alternative du même résultat) et de l'auteur \cite{marco2}.

\subsection{Notations} Soit $p$ un nombre premier. On fixe une clôture algébrique $\overline{\mathbb{Q}}_p$ de $\Q$ et une extension finie $F$ de $\Q$ contenue dans $\overline{\mathbb{Q}}_p$. On désignera toujours par $E$ une extension finie de $\Q$ qui vérifie:
\[
|S| = [F:\Q],
\]
où $S \ugu \mathrm{Hom}_{alg}(F,E)$. 

En général, si $L$ désigne $F$ ou $E$, on note $\mathcal{O}_L$ son anneau d'entiers, $\varpi_L$ une uniformisante de $\mathcal{O}_L$ et $k_L = \mathcal{O}_L/(\varpi_L)$ son corps résiduel. On note $f = [k_F : \mathbb{F}_p]$, $q= p^f$ et $e$ l'indice de ramification de $F$ sur $\Q$, de sorte que $[F:\Q] = ef$ et $k_F \simeq \mathbb{F}_q$.  

La valuation $p$-adique $val_F$ sur $\overline{\mathbb{Q}}_p$ est normalisée par $val_F(p) = [F:\Q]$ et on pose $|x| = p^{- val_F(x)}$  si $x \in \overline{\mathbb{Q}}_p$.

Si $a \in F$ et $n\in \mathbb{Z}$ on note $D(a,n) = a + \varpi_F^n\OF$, le disque de centre $a$ et de rayon $q^{-n}$.  

Soit $S'$ un sous-ensemble de $S$. Si $\underline{n}_{S'} = (n_{\sigma})_{\sigma \in S'}$,  $\underline{m}_{S'} = (m_{\sigma})_{\sigma \in S'}$ sont des $|S'|$-uplets d'entiers positifs ou nuls posons: 
\begin{itemize}
\item[(i)] $\underline{n}_{S'}! = \prod_{\sigma \in S'}n_{\sigma}!$;
\item[(ii)] $|\underline{n}_{S'}| = \sum_{\sigma \in S'}n_{\sigma}$;
\item[(iii)] $\underline{n}_{S'}-{\underline{m}_{S'}} = (n_{\sigma}-m_{\sigma})_{\sigma \in S'}$; 
\item[(iv)] $\underline{n}_{S'}\leqslant \underline{m}_{S'}$ si $n_{\sigma}\leq m_{\sigma}$ pour tout $\sigma \in S'$; 
\item[(v)] $\binom{\underline{n}_{S'}}{\underline{m}_{S'}}  \frac{\underline{n}_{S'}!}{\underline{m}_{S'}!(\underline{n}_{S'}-\underline{m}_{S'})!}$.
\end{itemize} 
Si $\underline{n}_{S'} = (n_{\sigma})_{\sigma \in S'} \in \mathbb{Z}_{\geq 0}^{|S'|}$ et $z \in \OF$ on pose $z^{\underline{n}_{S'}} = \prod_{\sigma \in S'}\sigma(z)^{n_{\sigma}}$. 

Pour alléger l'écriture, nous notons $\underline{n}$ un $|S|$-uplet d'entiers positifs ou nuls au lieu de $\underline{n}_S$.

Si $V$ est un $E$-espace vectoriel topologique, on note $V^{\vee}$ son dual topologique.

\subsection{Énoncé des résultats}

Pour énoncer le résultat principal il nous faut introduire un certain nombre de constructions. Soit $J$ une partie de $S$, $\underline{d}_{S\setminus J}$ un $|S\backslash J|$-uplet d'entiers positifs ou nuls. Posons:
\[
J' = J \coprod \{\sigma \in S\setminus J, d_{\sigma}+1 > -\val_{\Q}(\chi_1(p))\}.
\]
Soient $\chi_1, \chi_2$ deux caractères multiplicatifs localement $J$-analytiques de $F^{\times}$ dans $E^{\times}$. Notons $\chi_1\otimes \chi_2$ le caractère de $T$ défini par: 
\[
(\chi_1\otimes \chi_2)(\left[\begin{smallmatrix} {a} & {0} \cr {0} & {d} \end{smallmatrix}\right]) = \chi_1(a) \chi_2(d),
\]
où $T$ désigne le tore déployé constitué par les matrices diagonales de $G$. Par inflation on en déduit une représentation localement $J$-analytique de $P$. Notons:
\begin{itemize}
\item[$\bullet$] $\big(\Ind_{P}^G \chi_1\otimes \chi_2 \big)^{J-{an}}$ l'induite parabolique localement $J$-analytique, c'est-à-dire l'espace des fonctions localement $J$-analytiques $f$ sur $G$ à valeurs dans $E$ telles que $f(bg)=  (\chi_1\otimes\chi_2)(b)f(g)$ (l'action de $G$ étant la translation usuelle à droite sur les fonctions);
\item[$\bullet$]   $(\mathrm{Sym}^{d_{\sigma}}E^2)^{\sigma}$, pour $\sigma \in S$ et $d_{\sigma} \in \mathbb{Z}_{\geq 0}$, la représentation algébrique irréductible de $\mathrm{GL}_2 \otimes_{F,\sigma} E$ dont le plus haut poids est $\chi_{\sigma}\colon \mathrm{diag}(x_1,x_2) \mapsto \sigma(x_2)^{d_{\sigma}}$ vis-à-vis du sous-groupe des matrices triangulaires supérieures.
\end{itemize}
Considérons la représentation localement $\Q$-analytique de $G$ suivante: 
\[
I(\chi,J,\underline{d}_{S\setminus J})
 = \Big(  \bigotimes_{\sigma \in S \setminus J} (\Sym^{d_{\sigma}} E^2)^{\sigma}  \Big) \otimes_E \Big(\Ind_{P}^{G} \chi_1\otimes \chi_2 \Big)^{J-{an}}.
\]
Une première observation est que $I(\chi,J,\underline{d}_{S\setminus J})$ définit un faisceau sur $\mathbf{P}^1(F)$ dont les sections globales sont les fonctions $f\colon F \to E$ qui vérifient les deux conditions suivantes:  
\begin{itemize}
\item[(i)] $f|_{\OF}$ est une fonction dans $\mathcal{F}(\OF,J,\underline{d}_{S\setminus J})$  (Définition \ref{definloc});
\item[(ii)] $\chi_2 \chi_1^{-1}(z) z^{\underline{d}_{S\setminus J}}   f(1/z)|_{\OF-\{0\}}$ se prolonge sur $\OF$ en une fonction dans $\mathcal{F}(\OF,J,\underline{d}_{S\setminus J})$.
\end{itemize} 
Par ailleurs, des formules explicites munissent ce faisceau d'une action continue de $G$. D'après la preuve de \cite[Proposition 1.21]{eme}, le complété unitaire universel de $I(\chi,J,\underline{d}_{S\setminus J})$ est le complété par rapport au sous-$\Oe[P]$-réseau engendré par les vecteurs:
\[
\mathbf{1}_{\OF}(z)z^{\underline{n}_{S\setminus J}}z^{\underline{m}_{J}}, \quad \mathbf{1}_{F-\OF}(z)\chi_2\chi_1^{-1}(z) z^{\underline{d}_{S\setminus J}-\underline{n}_{S\setminus J}}z^{-\underline{m}_{J}}
\]
pour tout $\underline{0}\leqslant\underline{n}_{S\setminus J} \leqslant \underline{d}_{S\setminus J}$ et $\underline{m}_{J} \in \mathbb{Z}_{\geq 0}^{|J|}$. Notons $I(\chi,J,\underline{d}_{S\setminus J})^{\bigwedge}$ le complété de $I(\chi,J,\underline{d}_{S\setminus J})$ par rapport à ce réseau. 

Le but de cet article est celui de donner une description explicite de l'espace $I(\chi,J,\underline{d}_{S\setminus J})^{\bigwedge}$. Dans un premier temps nous démontrons deux résultats qui ajoutent des conditions supplémentaires aux données initiales et qui permettent de ne pas considérer des cas pathologiques ou bien de simplifier le problème. Le premier ingrédient donne deux conditions nécessaires pour que $I(\chi,J,\underline{d}_{S\setminus J})^{\bigwedge}$ soit non nul.

\begin{prop} \label{nullita11}
Le deux conditions suivantes sont nécessaires pour que   $I(\chi,J,\underline{d}_{S\setminus J} )^{\bigwedge}$ soit non nul: 
\begin{itemize}
\item[(i)] Le caractère central de $I(\chi,J,\underline{d}_{S\setminus J} )$ est intègre; 
\item[(ii)] On a l'inégalité $\mathrm{val}_{\Q}(\chi_2(p)) + |\underline{d}_{S\setminus J}| \geq 0$.  
\end{itemize}
\end{prop}
Mentionnons qu'il s'agit d'un résultat bien connu pour $F=\Q$ (\cite[Lemma 2.1]{eme}) et en dehors de $\Q$ dans le cas localement algébrique, c'est-à-dire $J=\emptyset$ (\cite[Lemme 7.9]{pas}). En particulier, si les conditions de la Proposition \ref{nullita11} sont satisfaites on déduit que $r \geq 0$ où   $r \ugu -\mathrm{val}_{\Q}(\chi_1(p))$.

Notons $\chi_1' = \chi_1$, $\chi_2' = \chi_2\, \prod_{\sigma \in {J'\setminus J}} \sigma^{d_{\sigma}}$ 
et remarquons que l'on a une immersion fermée $G$-équivariante:
\begin{align}\label{immersionech}
I(\chi,J,\underline{d}_{S\setminus J}) \into I(\chi',J',\underline{d}_{S\setminus J'}). 
\end{align}
Un autre ingrédient important est la proposition suivante, essentiellement démontre par Breuil en faisant recours aux techniques dévéloppées par Amice-Vélu et Vishik, qui  donne des indications concernant la structure de $I(\chi,J,\underline{d}_{S\setminus J} )^{\bigwedge}$, ou plus précisement ses vecteurs localement $\Q$-analytiques. 

\begin{prop}\label{dsigma12}
Supposons que les conditions de la Proposition \ref{nullita11} soient satisfaites. Alors les conditions suivantes sont équivalentes:
\begin{itemize}
\item[(i)] Toute application continue, $E$-linéaire et $G$-équivariante $I(\chi,J,\underline{d}_{S\setminus J}) \to B$, où $B$ est un $G$-Banach unitaire, s'étend de manière unique en une application continue, $E$-linéaire et $G$-équivariante $I(\chi',J',\underline{d}_{S\setminus J'}) \to B$. 
\item[(ii)] L'application canonique $I(\chi,J,\underline{d}_{S\setminus J}) \to I(\chi,J,\underline{d}_{S\setminus J} )^{\bigwedge}$ s'étend de manière unique en une application continue, $E$-linéaire et $G$-équivariante $I(\chi',J',\underline{d}_{S\setminus J'}) \to I(\chi,J,\underline{d}_{S\setminus J} )^{\bigwedge}$.  
\item[(iii)] L'application \eqref{immersionech} induit un isomorphisme de $G$-Banach unitaires: 
\[
I(\chi,J,\underline{d}_{S\setminus J} )^{\bigwedge} \stackrel{\sim}{\longrightarrow} I(\chi',J',\underline{d}_{S\setminus J'})^{\bigwedge}
\]
\end{itemize} 
\end{prop}

Donc, d'après la Proposition \ref{dsigma12} (iii) on est réduit à considérer  $I(\chi',J',\underline{d}_{S\setminus J'} )^{\bigwedge}$. Par un calcul analogue à celui dans la preuve de \cite[Théorème 4.3.1]{bb} on trouve qu'une boule ouverte (de centre $0$) du Banach dual de $I(\chi',J',\underline{d}_{S\setminus J'})^{\bigwedge}$ s'identifie aux distributions $\mu$ dans le dual fort de $I(\chi',J',\underline{d}_{S\setminus J'})$ telles que pour tout $n \in \mathbb{Z}$, tout $a \in F$, tout $\underline{0}\leqslant \underline{n}_{S\setminus J'}\leqslant \underline{d}_{S\setminus J'}$ et tout $\underline{m}_{J'} \in \mathbb{Z}_{\geq 0}^{|J|}$ on a:
\begin{align}
&\Big|\int_{D(a,n)} (z-a)^{\underline{n}_{S\setminus J'}}  (z-a)^{\underline{m}_{J'}} \mu(z)\Big| \leq C_{\mu} q^{n (r-|\underline{n}_{S\setminus J'}|-|\underline{m}_{J'}|)}  \label{puno13} \\
&\Big|\int_{F\setminus D(a,n+1)} \chi_2\chi_1^{-1}(z-a) (z-a)^{\underline{d}_{S\setminus J'}- \underline{n}_{S\setminus J'}} (z-a)^{-\underline{m}_{J'}} \mu(z)\Big| \leq C_{\mu} q^{n(|\underline{n}_{S\setminus J'}|+|\underline{m}_{J'}|-r)},  \label{pdue13}
\end{align}
où $C_{\mu} \in \mathbb{R}_{\geq 0}$.

D'autre part, un étude fine du dual fort de l'espace de Banach des fonctions de classe $C^r$ sur $\OF$ ou plus précisement de son sous-espace fermé $C^r(\OF,J',\underline{d}_{S\setminus J'})$ (§\ref{chiusi}) nous fournit une condition nécessaire et suffisante pour qu'une forme linéaire sur $\mathcal{F}^N(\OF,J,\underline{d}_{S\setminus J})$ (voir §\ref{ana} pour une définition de cet espace) s'étende en une distribution sur  $C^r(\OF,J',\underline{d}_{S\setminus J'})$ (Théorème \ref{velu}). Pour $F = \Q$ il s'agit d'un résultat bien connu et d\^u à Amice-Vélu et Vishik (\cite{amivel},  \cite{vis}). Plus précisément: 

\begin{theo}\label{velu7}  (i) Soit $\mu \in C^r(\OF,J', \underline{d}_{ S \setminus J'})^{\vee}$. Il existe une constante $C_{\mu} \in \mathbb{R}_{\geq 0}$ telle que pour tout $a \in \OF$, tout $n \in \mathbb{Z}_{\geq 0}$, tout $\underline{0} \leqslant \underline{n}_{S\setminus J'} \leqslant \underline{d}_{S\setminus J'}$ et tout $\underline{m}_{J'} \in \mathbb{Z}_{\geq 0}^{|J'|}$ on ait: 
\begin{align*}
\Big|\int_{D(a,n)} (z-a)^{\underline{n}_{S\setminus J'}} (z-a)^{\underline{m}_{J'}} \mu(z)\Big| \leq C_{\mu} \, q^{n(r-|\underline{n}_{S\setminus J'}| - |\underline{m}_{J'}|)}. 
\end{align*}

(ii) Soit $N$ un entier tel que $N\geq [r]$ et $\mu$ une forme linéaire sur  $\mathcal{F}^N(\OF,J, \underline{d}_{ S \setminus J})$. Supposons qu'il existe une constante $C_{\mu} \in \mathbb{R}_{\geq 0}$ telle que pour tout $a \in \OF$, tout $n \in \mathbb{Z}_{\geq 0}$, tout $\underline{0} \leqslant \underline{n}_{S\setminus J} \leqslant \underline{d}_{S\setminus J}$ et tout $\underline{m}_{J} \in \mathbb{Z}_{\geq 0}^{|J|}$ tels que $|\underline{n}_{S\setminus J}|+|\underline{m}_{J}|\leq N$, on ait:    
\begin{align*}
\Big|\int_{D(a,n)} (z-a)^{\underline{n}_{S\setminus J}} (z-a)^{\underline{m}_{J}} \mu(z)\Big| \leq C_{\mu} \, q^{n(r-|\underline{n}_{S\setminus J}| - |\underline{m}_{J}|)}.
\end{align*}
Alors $\mu$ se prolonge de manière unique en une distribution sur $C^r(\OF,J', \underline{d}_{ S \setminus J'})$. 
\end{theo}

On est alors amené à considérer l'espace $B(\chi',J',\underline{d}_{S\setminus J'})$ des fonctions $f$ de $F$ dans $E$ qui vérifient les deux conditions suivantes: 
\begin{itemize}
\item[(i)] $f|_{\OF}$ est une fonction dans $C^{r}(\OF,J',\underline{d}_{S\setminus J'})$;
\item[(ii)] $\chi_2' {\chi_1'}^{-1}(z) z^{\underline{d}_{S\setminus J'}} f(1/z)|_{\OF-\{0\}}$ se prolonge sur $\OF$ en une fonction dans $C^{r}(\OF,J',\underline{d}_{S\setminus J'})$,
\end{itemize}
qui est un espace de Banach $p$-adique naturellement muni d'une action continue de $G$.

Un examen approfondi, qui utilise de manière cruciale le Théorème \ref{velu7}, montre que les conditions \eqref{puno13} et \eqref{pdue13} sélectionnent exactement les formes linéaires dans $B(\chi',J',\underline{d}_{S\setminus J'})^{\vee}$  annulant les fonctions d'un sous-espace $L(\chi',J',\underline{d}_{S\setminus J'})$ de $B(\chi',J',\underline{d}_{S\setminus J'})$ (voir §\ref{repban} pour une définition de cet espace).

Le résultat principal de cet article, qui généralise le \cite[Théorème 4.3.1]{bb} pour $F=\Q$, est alors le suivant.

\begin{theo}\label{teopr}
Il existe un isomorphisme $G$-équivariant d'espaces de Banach $p$-adiques: 
\[
I(\chi,J,\underline{d}_{S\setminus J})^{\bigwedge} \stackrel{\sim}{\longrightarrow} B(\chi,J',\underline{d}_{S\setminus J'})/L(\chi,J',\underline{d}_{S\setminus J'}).
\]
\end{theo}

\subsection{Plan de l'article} Dans la Section $2$ nous rappelons quelques généralités d'analyse fonctionnelle $p$-adique et la notion de complété unitaire universel introduite dans \cite{eme}.   La Section $3$ est constituée de quelques rappels sur les espaces des fonctions de classe $C^r$ et ses duaux. Nous introduisons dans la Section $4$ les représentations localement $\Q$-analytiques $I(\chi,J,\underline{d}_{S\setminus J})$ qui font l'objet de notre étude et ensuite nous construisons la représentation de Banach $\Pi(\chi,J,\underline{d}_{S\setminus J})$. Dans la section $5$ nous donnons deux conditions nécessaires pour que le complété unitaire universel de $I(\chi,J,\underline{d}_{S\setminus J})$ soit non nul et ensuite nous commençons l'étude des espaces duaux $(I(\chi,J,\underline{d}_{S\setminus J})^{\bigwedge})^{\vee}$ et $\Pi(\chi,J,\underline{d}_{S\setminus J})^{\vee}$. Dans la section $6$ qui est le c\oe ur de cet article, nous démontrons le Théorème \ref{teopr} et on termine avec un exemple explicite.  




\addtocontents{toc}{\protect\setcounter{tocdepth}{2}}


\section{Préliminaires}

\subsection{Rappels d'analyse fonctionnelle non archimédienne}\label{funzionale}

Dans ce paragraphe on donne divers notions préliminaires d'analyse fonctionnelle non archimédienne dont on se servira par la suite. Nous renvoyons à \cite{sch} pour plus de détails. 

Un $E$-espace vectoriel topologique $V$ est dit \textit{localement convexe} si une base de voisinages de l'origine peut être définie par une famille de sous-$\Oe$-modules de $V$. Ou de manière équivalente si la topologie peut être définie par une famille de semi-normes non archimédiennes \cite[Proposition 4.3, Proposition 4.4]{sch}. 

Soit $V$ un $E$-espace vectoriel localement convexe. Un \textit{réseau} $\mathcal{L}$ de $V$ est un sous-$\Oe$-module de $V$ tel que pour tout $v \in V$ il existe un élément non  nul $a \in E^{\times}$ tel que $av \in \mathcal{L}$.  En particulier, on observe que tous les sous-$\Oe$-modules ouverts de $V$ sont des réseaux. Si l'on se donne deux réseaux $\mathcal{L}_1$ et $\mathcal{L}_2$ de $V$, on dit qu'ils sont \textit{commensurables} s'il existe $a \in E^{\times}$ tel que $a \mathcal{L}_1 \subseteq \mathcal{L}_2 \subseteq a^{-1} \mathcal{L}_2$. La commensurabilité définit une relation d'équivalence sur l'ensemble $\mathcal{L}(V)$ des réseaux ouverts. 

On dit qu'un réseau $\mathcal{L}$ de $V$ est \textit{séparé}  si $\bigcap_{n \in \mathbb{N}} \varpi_E^n \mathcal{L} = 0$ ou, de manière équivalente, s'il ne contient pas de $E$-droite.

On dit que $V$ est \textit{tonnelé} si tout réseau fermé dans $V$ est ouvert.

On dit que $V$ est de \textit{Fréchet} s'il est complet et métrisable ou, de manière équivalente, s'il est complet, Hausdorff, et sa topologie peut être définie par une famille dénombrable de semi-normes. En particulier, si sa topologie peut être définie par une unique norme, on dit que $V$ est un espace de \textit{Banach}. 

Si $V$ est un espace de Banach sur $E$ alors un réseau ouvert est séparé si et seulement s'il est borné. En outre, si $\mathcal{L}$ est un réseau ouvert et séparé de $V$ alors la jauge de $\mathcal{L}$ définie par:
\[
\forall v \in V, \quad   \|v\|_{\mathcal{L}} = \inf_{v \in a\mathcal{L}} |a|
\]
est une norme et la topologie sur $V$ définie par $\|\cdot \|_{\mathcal{L}}$ coïncide avec celle initiale \cite[Corollaire 4.12]{sch}.

On dit que $V$ est de \textit{type compact} s'il existe un isomorphisme de $E$-espaces vectoriels topologiques:
\[
V \longrightarrow \varinjlim_{n} V_n
\]
où $\{V_n \}_{n\geq 1}$ est un système inductif d'espaces de Banach sur $E$, tel que les morphismes de transition sont injectifs et compacts.

Un sous-ensemble $B \subseteq V$ est dit \textit{borné} si pour tout réseau $\mathcal{L} \subseteq V$ il existe $a \in E$ tel que $B \subseteq a\mathcal{L}$. 

Soit $W$ un $E$-espace vectoriel localement convexe. On note $\Hom_E(V,W)$ l'espace des fonctions $E$-linéaires et continues sur $V$ à valeurs dans $W$. On fixe un sous-ensemble borné $B \subseteq V$. Si $p$ est une semi-norme continue sur $W$ alors la formule: 
\[
p_B(f) = \sup_{v \in B} p(f(v))
\] 
définit une semi-norme sur $\Hom_E(V,W)$. Soit $\mathcal{B}$ une famille de sous-ensembles bornés de $V$. La topologie localement convexe sur $\Hom_E(V,W)$ définie par la famille de semi-norme $\{p_B : B \in \mathcal{B}, p \ \mbox{semi-norme continue sur} \ W \}$ est appelée $\mathcal{B}$-topologie. Si $\mathcal{B}$ est la famille de tous les singletons alors la $\mathcal{B}$-topologie correspondante est aussi dite \textit{topologie faible}. Si $\mathcal{B}$ est la famille de tous les sous-ensembles bornés de $V$ alors la $\mathcal{B}$-topologie correspondante est dite \textit{topologie forte}.

 \subsection{Complétés unitaires universels}

Soit $G$ le groupe des $\Q$-points d'un groupe algébrique linéaire réductif connexe défini sur $\Q$. La notion de complété unitaire universel pour un espace vectoriel localement convexe  muni d'une action continue de $G$ a été formalisée par Emerton dans \cite[§1]{eme}, après que des exemples de complétés unitaires universels aient été construits par Breuil (\cite{breuila, breuilb}) et Berger-Breuil (\cite{bb}). Dans ce paragraphe nous rappelons le contexte dans lequel cette notion s'insère ainsi qu'une condition nécessaire et suffisante pour qu'un tel complété existe.

\begin{defi}[\cite{sch3, breuilb}]\label{banachunit} Un $G$-Banach est un espace de Banach $B$ sur $E$ muni d'une action à gauche de $G$ telle que l'application $G \times B \to B$ décrivant cette action est continue. Un $G$-Banach $B$ est dit unitaire si pour un choix de norme $\|\cdot \|$ définissant la topologie de $B$, on a $\|gv \| = \|v \|$ pour tout $g \in G$ et tout $v \in B$.
\end{defi} 

\begin{rema}{\rm
Si le groupe $G$ est compact alors tout $G$-Banach est unitaire. Ceci n'est pas vrai si $G$ n'est pas compact.}
\end{rema}

Soit $V$ un $E$-espace vectoriel localement convexe et muni d'une action continue de $G$. Un complété unitaire universel de $V$ est un $G$-Banach unitaire $B$ qui satisfait une certaine propriété universelle. Plus précisément:

\begin{defi}[\cite{eme}, Définition 1.1]\label{comple}
Avec les notations précédentes, un complété unitaire universel de $V$ est la donnée d'un $G$-Banach unitaire $B$ et d'une application $E$-linéaire, continue et $G$-équivariante $\iota$ de $V$ sur $B$ telle que toute application $E$-linéaire, continue et $G$-équivariante $V \to W$, où $W$ est un $G$-Banach unitaire, se factorise de façon unique à travers $\iota$.
\end{defi}

\begin{rema}\label{ossuti} 
{\rm Si un complété unitaire universel $(B,\iota)$ de $V$ existe, alors il est unique à isomorphisme près. Par ailleurs, l'adhérence dans $B$ de l'image de $V$ à travers $\iota$ vérifie la propriété universelle énoncée dans la Définition \ref{comple}. On en déduit que l'application $\iota$ est d'image dense. }
\end{rema}

Le lemme suivant fournit une condition nécessaire et suffisante pour qu'un tel complété unitaire universel existe.
\begin{lemm}[\cite{eme}, Lemme 1.3]\label{lemmaeme}
La $G$-représentation $V$ admet un complété unitaire universel si et seulement si l'ensemble des classes de commensurabilité des réseaux  ouverts et stables sous l'action de $G$ dans $V$, partiellement ordonné par l'inclusion, possède un élément minimal. 
\end{lemm}

\section{Rappels sur les fonctions de classe $\Cr$ sur $\OF$}

Soit $r \in \mathbb{R}_{\geq 0}$. Dans \cite{marco} nous avons introduit une nouvelle notion de fonction de classe $C^r$ sur $\OF$, qui s'appuie principalement sur les travaux d'Amice, Amice-Velù, Vishik, Van der Put et Colmez (\cite{ami}, \cite{amivel}, \cite{vis}, \cite{vander}, \cite{colmez}). Dans cette section nous allons rappeler un certain nombre de constructions et de résultats concernant l'espace des fonctions de classe $C^r$ sur $\OF$. Nous renvoyons à \cite{marco} pour plus de détails et à \cite{enno1} et \cite{enno2} pour d'autres possibles  notions.  

\subsection{Définition et compléments}\label{classe}

Soit $r \in \mathbb{R}_{\geq 0}$. Notons $[r]$ sa partie entière. Si $n \in \mathbb{Z}_{\geq 0}$ et $* \in \{<,\leq,>, \geq, = \}$ notons: 
\[
I_{*n} = \Big\{\underline{i} \in \mathbb{Z}_{\geq 0}^{|S|}, \, \sum_{\sigma \in S}i_{\sigma} * n  \Big\}.
\]

\begin{defi}\label{definizione}
On dit que $f\colon \OF \to E$ est de classe $\Cr$ sur $\OF$ s'il existe des fonctions bornées $D_{\underline{i}}f \colon \OF \to E$, pour $\underline{i} \in I_{\leq [r]}$, telles que, si l'on définit $\varepsilon_{f,[r]}\colon \OF \times \OF \to E$ par: 
\begin{align*}
\varepsilon_{f,[r]}(x,y) = f(x+y) -
 \sum_{\underline{i} \in I_{\leq [r]}} D_{\underline{i}}f(x) \frac{y^{\underline{i}}}{\underline{i}!}
\end{align*}
et pour tout $h \in \mathbb{Z}_{\geq 0}$
\begin{align*}
C_{f,r}(h) = \sup_{x \in \OF, y \in \varpi_F^h \OF} |\varepsilon_{f,[r]}(x,y)| \, q^{rh}
\end{align*}
alors $C_{f,r}(h)$ tend vers $0$ quand $h$ tend vers $+\infty$. 
\end{defi}

Si $f$ est une fonction de classe $C^r$ sur $\OF$ alors il existe une unique famille de fonctions 
\[
\big\{ D_{\underline{i}}f \colon \OF \to E, \ \underline{i} \in I_{\leq [r]} \big\}
\]
qui vérifie la Définition \ref{definizione} (\cite[Lemme 2.4]{marco}). Notons $C^r(\OF,E)$ l'ensemble des fonctions de $\OF$ dans $E$ qui sont de classe $C^r$ et munissons-le de la norme $\|\cdot \|_{C^r}$ définie par:
\[
\|f \|_{C^r} = \sup \Big( \sup_{\underline{i} \in I_{\leq [r]}} \sup_{x \in \OF} \Bigl|\frac{D_{\underline{i}}f(x)}{\underline{i}!}\Bigr|, \sup_{x,y \in \OF} \frac{| \varepsilon_{f,[r]}(x,y)|}{|y|^{r}}   \Big)
\]
ce qui en fait un espace de Banach sur $E$. Plus précisément l'espace $C^r(\OF,E)$ est une \textit{$E$-algèbre de Banach} (\cite[Lemme 2.9]{marco}), c'est-à-dire une $E$-algèbre normée telle que l'espace vectoriel normé sous-jacent soit un espace de Banach.

Montrons le résultat suivant dont on se servira par la suite. 

\begin{lemm}\label{corollarioutile}
Soit $n \in \mathbb{Z}_{\geq 0}$ et $f$ une fonction de classe $C^r$ sur $\OF$. Notons $g$ la fonction de $\OF$ dans $E$ définie par:
\[
z \mapsto \mathbf{1}_{D(0,n)}(z) f\Big(\frac{z}{\varpi_F^n} \Big).
\]
Alors $g \in C^r(\OF,E)$ et $\|g \|_{C^r} \leq q^{nr} \|f \|_{C^r}$.
\begin{proof}
Posons pour tout $\underline{i} \in I_{\leq [r]}$:   
\begin{align}\label{definizionederiv}
\forall z \in \OF,\quad D_{\underline{i}}g(z) = \Big(\frac{1}{\varpi_F^n}  \Big)^{\underline{i}} \mathbf{1}_{\varpi_F^n \OF}(z)D_{\underline{i}}f \Big(\frac{z}{\varpi_F^n} \Big)
\end{align} 
et
\[
\forall x,y \in \OF, \quad \varepsilon_{g,[r]}(x,y) = \mathbf{1}_{D(0,n)}(x+y) f\Big(\frac{x+y}{\uni^n} \Big) - \sum_{\underline{i}\in I_{\leq [r]}} \mathbf{1}_{D(0,n)}(x) D_{\underline{i}}f \Big(\frac{x}{\uni^n} \Big) \Big(\frac{y}{\uni^n} \Big)^{\underline{i}}. 
\]
On voit immédiatement que l'on a
\[
\forall h\geq n, \quad \sup_{x\in \OF, y\in \uni^h\OF}|\varepsilon_{g,[r]}(x,y)| \leq \sup_{x\in \OF, y\in \uni^h\OF}|\varepsilon_{f,[r]}(x,y)|
\]
ce qui implique que $g$ est de classe $C^r$. Il nous reste à montrer l'inégalité  sur la norme. Par \eqref{definizionederiv} on a:
\begin{align}\label{inegcorut}
\forall \underline{i} \in I_{\leq [r]}, \quad \sup_{z \in \OF} \Big| \frac{D_{\underline{i}}g(z)}{\underline{i}!} \Big| \leq \Big|\Big(\frac{1}{\varpi_F^n}\Big)^{\underline{i}}\Big| \sup_{z \in \OF} \Big| \frac{D_{\underline{i}}f(z)}{\underline{i}!} \Big| \leq q^{n|\underline{i}|} \|f \|_{C^r} \leq q^{nr} \|f \|_{C^r}.
\end{align}
Par ailleurs:
\begin{itemize}
\item[$\bullet$] Dans le cas $x,y \in \varpi_F^n \OF$ on a:
\[
 \frac{| \varepsilon_{g,[r]}(x,y)|}{|y|^{r}} \leq \frac{| \varepsilon_{f,[r]}(x,y)|}{|y|^{r}} \leq q^{nr} \|f \|_{C^r}.
\]
\item[$\bullet$] Dans le cas $x \in \varpi_F^n \OF$, $y \notin \varpi_F^n \OF$ on a:
\begin{align*}
 \frac{| \varepsilon_{g,[r]}(x,y)|}{|y|^{r}} = \frac{\big|\sum_{\underline{i}\in I_{\leq [r]}}D_{\underline{i}}g(x) \frac{y^{\underline{i}}}{\underline{i}!}\big|}{|y|^{r}}
 &\leq \sup_{\underline{i}\in I_{\leq [r]}} \Big|\frac{D_{\underline{i}}g(x)}{\underline{i}!}  \Big| |y|^{|\underline{i}|-r} \\
 &\leq \sup_{\underline{i}\in I_{\leq [r]}} \sup_{x \in \OF} \Big|\frac{D_{\underline{i}}f(x)}{\underline{i}!}  \Big| |\varpi_F^n|^{-|\underline{i}|}|y|^{|\underline{i}|-r} \\
 &\leq q^{nr} \|f\|_{C^r}.
\end{align*}
\item[$\bullet$] Dans le cas $x \notin \varpi_F^n \OF$, $x+y \notin \varpi_F^n \OF$ on a $\varepsilon_{g,[r]}(x,y) = 0$.
\item[$\bullet$] Dans le cas $x \notin \varpi_F^n \OF$, $x+y \in \varpi_F^n \OF$ on a:
\begin{align*}
\frac{| \varepsilon_{g,[r]}(x,y)|}{|y|^{r}} = \frac{\Big|f\Big(\frac{x}{\varpi_F^n}+ \frac{y}{\varpi_F^n} \Big)  \Big|}{|y|^{r}}
 \leq  q^{nr} \|f\|_{C^r}.
\end{align*}
\end{itemize}
Ceci permet, en revenant à la définition de $\|\cdot \|_{C^r}$, de conclure.
\end{proof}
\end{lemm}

\subsubsection{\textbf{Composition de fonctions}}
Soit $f$ une fonction de $\OF$ dans $E$ de classe $C^r$ et $h$ une fonction de $\OF$ dans $\OF$. Dans ce paragraphe nous rappelons (\cite[§2.2.1]{marco}) une condition suffisante sur $h$ pour que $f\circ h\colon \OF \to E$ soit de classe $C^r$. Pour cela, nous avons besoin d'introduire la définition suivante.
\begin{defi}\label{definizione22}
Soit $r \in \mathbb{R}_{\geq 0}$. On dit que $h\colon \OF \to F$ est de classe $C^{r,id}$ sur $\OF$ s'il existe des fonctions bornées $h^{(i)} \colon \OF \to F$, pour $0\leq i \leq [r]$, telles que, si l'on définit $\varepsilon_{h,[r]}\colon \OF \times \OF \to F$ par: 
\begin{align*}
\varepsilon_{h,[r]}(x,y) = f(x+y) -
 \sum_{i = 0}^{[r]} h^{(i)}(x) \frac{y^{i}}{i!}
\end{align*}
et pour tout $k \in \mathbb{Z}_{\geq 0}$
\begin{align*}
C_{h,r}(k) = \sup_{x \in \OF, y \in \varpi_F^k \OF} |\varepsilon_{h,[r]}(x,y)| \, q^{rk}
\end{align*}
alors $C_{h,r}(k)$ tend vers $0$ quand $k$ tend vers $+\infty$. 
\end{defi}

Notons $C^{r,id}(\OF,F)$ l'ensemble des fonctions de $\OF$ dans $F$ qui sont de classe $C^{r,id}$. On munit $C^{r,id}(\OF,F)$ de la norme $\|\cdot \|_{C^{r,id}}$ définie par:
\[
\|h \|_{C^{r,id}} = \sup \Big( \sup_{0 \leq i \leq [r]} \sup_{x \in \OF} \Bigl|\frac{h^{(i)}(x)}{i!}\Bigr|, \sup_{x,y \in \OF} \frac{| \varepsilon_{h,[r]}(x,y)|}{|y|^{r}}   \Big)
\]
ce qui en fait un espace de Banach sur $F$.

\begin{prop}\label{composizione}
Soit $r \in \mathbb{R}_{\geq 0}$. Si $h\colon \OF \to \OF$ est une fonction de classe $C^{r,id}$ alors
\begin{itemize}
\item[(i)] $\forall f \in C^r(\OF,E)$,  $f \circ h \in C^r(\OF,E)$;
\item[(ii)] l'application de $C^r(\OF,E)$ dans $C^r(\OF,E)$ définie par $f \mapsto f \circ h$ est continue.
\end{itemize}
\begin{proof}
Voir \cite[Proposition 2.12]{marco}.
\end{proof}
\end{prop}

\subsubsection{\textbf{Construction de sous-espaces fermés}}\label{chiusi}
Soit $r \in \mathbb{R}_{\geq 0}$, $J \subseteq S$ et $d_{\sigma} \in \mathbb{Z}_{\geq 0}$ pour $\sigma \in S\backslash J$. Nous allons définir un sous-espace fermé de $C^r(\OF,E)$ qui dépend de $J$ et de $\underline{d}_{S\setminus J}$ et qui va jouer un rôle important dans la suite. 

Posons:
\[
J' = J \coprod \{\sigma \in S\backslash J,\, d_{\sigma}+1 > r\}
\] 
et désignons par $e_{\sigma}$ le vecteur de $\mathbb{Z}_{\geq 0}^{|S|}$ ayant toutes ses composantes nulles sauf celle d'indice $\sigma$ qui est égal à $1$. Notons pour tout $f \in C^r(\OF,E)$:   
\[
\forall \sigma \in S,\ 0\leq i \leq [r], \quad \frac{\partial^i }{\partial z_{\sigma}^i}f = D_{i e_{\sigma}}f.
\] 

\begin{defi}\label{sottosp}
On note $C^r(\OF,J',\underline{d}_{S\setminus J'})$ le sous-$E$-espace vectoriel des fonctions $f$ de classe $C^r$ sur $\OF$ telles que: 
\[
\forall \sigma \in S\backslash J', \quad \frac{\partial^{d_{\sigma}+1} }{\partial z_{\sigma}^{d_{\sigma}+1}}f = 0. 
\] 
\end{defi} 

D'après \cite[Corollaire 2.8]{marco} l'opérateur $D_{\underline{i}}$ est continu pour tout $\underline{i} \in I_{\leq [r]}$ ce qui implique que l'espace $C^r(\OF,J',\underline{d}_{S\setminus J'})$ est bien un sous-espace fermé de $C^r(\OF,E)$.  On le munit de la topologie induite par celle de $C^r(\OF,E)$ qui en fait un espace de Banach sur $E$.

\subsection{Fonctions localement analytiques et fonctions de classe $C^r$ }\label{ana}

Soit $U\subseteq \OF$ un sous-ensemble ouvert, $J \subseteq S$ et $d_{\sigma} \in \mathbb{Z}_{\geq 0}$ pour $\sigma \in S\backslash J$. Pour $a \in U$ et $n \in \mathbb{Z}_{\geq 0}$ tels que $D(a,n) \subseteq U$, on note $\mathcal{O}(D(a,n),J, \underline{d}_{S\setminus J})$ le $E$-espace vectoriel des fonctions $f\colon D(a,n) \to E$ telles que 
\begin{align*}
f(z) = \sum_{\substack{
\underline{m} = (m_{\sigma})_{\sigma \in S} \in \mathbb{Z}_{\geq 0}^{|S|} \\
m_{\sigma} \leq d_{\sigma} \ \mbox{si} \ \sigma \in S\setminus J
}} a_{\underline{m}}(a) (z-a)^{\underline{m}}
\end{align*}
avec $a_{\underline{m}}(a) \in E$ et $|a_{\underline{m}}(a)| q^{-n(|\underline{m}|)} \to 0$ quand $|\underline{m}| \to +\infty$. Muni de la topologie définie par la norme  
\[
\| f \|_{a, n} = \sup_{\underline{m}} \left( |a_{\underline{m}}(a)| q^{-n (|\underline{m}|)} \right) 
\]
c'est un espace de Banach sur $E$. 

Par compacité de $U$ il existe $h_0 \in \mathbb{Z}_{\geq 0}$ tel que
\[
\forall a \in U, \forall h\geq h_0, \quad D(a,h) \subseteq U.
\]
Pour tout $h\geq h_0$ on note $\mathcal{F}_{h}(U,J, \underline{d}_{S\setminus J})$ le $E$-espace vectoriel des fonctions $f\colon U \to E$ telles  que:
\[
\forall a \in U, \quad f|_{D(a,h)} \in \mathcal{O}(D(a,h),J, \underline{d}_{S\setminus J}).
\]
On munit cet espace de la norme définie par:
\begin{align}\label{formulanorma}
\| f \|_{\mathcal{F}_h} = \sup_{a \, \mathrm{mod} \, \uni^h, a\in U} \|f|_{D(a,h)} \|_{a,h}
\end{align}
qui en fait un espace de Banach sur $E$. On voit immédiatement que cette définition ne dépend pas du choix des représentants. De plus (\cite[p. 107]{sch}) les inclusions 
\[
\mathcal{F}_{h}(U, J, \underline{d}_{S\setminus J}) \into \mathcal{F}_{h+1}(U, J, \underline{d}_{S\setminus J})
\] 
sont continues et compactes.

\begin{defi}\label{definloc}
On note $\mathcal{F}(U,J, \underline{d}_{S\setminus J})$ le $E$-espace vectoriel des fonctions $f\colon \OF \to E$ telles qu'il existe un entier $h$ tel que $h\geq h_0$ et
\[
f \in \mathcal{F}_{h}(U, J, \underline{d}_{S\setminus J}). 
\]
\end{defi}

On munit  l'espace $\mathcal{F}(U,J, \underline{d}_{S\setminus J})$ de la topologie de la limite inductive qui en fait un espace de type compact. Posons pour tout $N \in \mathbb{Z}_{\geq 0}$:
\begin{align*}
\mathcal{F}^{N}(\OF,S) & = \textstyle{\sum_{
\underline{d} \in I_{\leq N}}} \mathcal{F}(\OF,\emptyset, \underline{d}); \\
\mathcal{F}^{N}(\OF,J, \underline{d}_{S \setminus J}) &= \mathcal{F}^{N}(\OF, S) \cap \mathcal{F}(\OF,J, \underline{d}_{S \setminus J}). 
\end{align*} 
Notons que l'espace $\mathcal{F}^{N}(\OF,S)$ (resp. $\mathcal{F}^{N}(\OF,J, \underline{d}_{S \setminus J})$) est  un sous-$E$-espace vectoriel de $\mathcal{F}(\OF,S)$ (resp. $\mathcal{F}(\OF,J, \underline{d}_{S \setminus J})$) et rappelons les deux faits suivants:
\begin{itemize}
\item[$\bullet$] L'espace $\mathcal{F}(\OF,J, \underline{d}_{S \setminus J})$ s'injecte de façon continue dans $C^r(\OF,J',\underline{d}_{S\setminus J'})$ \cite[Corollaire 3.4]{marco};
\item[$\bullet$] Si $N$ est un entier tel que $N \geq [r]$, alors l'espace $\mathcal{F}^{N}(\OF,J, \underline{d}_{S \setminus J})$ est dense dans $C^r(\OF,J',\underline{d}_{S\setminus J'})$ \cite[Corollaire 3.16]{marco}.
\end{itemize}
En particulier, le deuxième point est conséquence du fait que l'on peut construire une base de Banach pour l'espace $C^r(\OF,J',\underline{d}_{S\setminus J'})$ qui est constituée de fonctions dans  $\mathcal{F}^{[r]}(\OF,J, \underline{d}_{S \setminus J})$.

\subsection{Distributions d'ordre $r$}\label{duali}
Conservons les notations du  §\ref{ana} et notons $\mathcal{F}^{N}(\OF,J, \underline{d}_{S \setminus J})^{\vee}$, pour tout $N\in \mathbb{Z}_{\geq 0}$,  l'ensemble des formes linéaires sur $\mathcal{F}^{N}(\OF,J, \underline{d}_{S \setminus J})$. Si $N$ est un entier tel que $N\geq [r]$ alors, d'après  \cite[Corollaire 3.16]{marco}, l'inclusion 
\[
\mathcal{F}^{N}(\OF,J,\underline{d}_{S \setminus J}) \subseteq C^r(\OF,J',\underline{d}_{S \setminus J'})
\] 
induit une injection
\[
C^r(\OF,J,\underline{d}_{S \setminus J})^{\vee} \into \mathcal{F}^{N}(\OF,J,\underline{d}_{S \setminus J})^{\vee}. 
\]
Dans cette section nous rappelons une condition nécessaire et suffisante pour qu'une forme linéaire  $\mu\colon \mathcal{F}^{N}(\OF,J,\underline{d}_{S \setminus J})\to E$ s'étende en une forme linéaire continue sur l'espace de Banach $C^r(\OF,J',\underline{d}_{S \setminus J'})$. Cela généralise un résultat d\^u à Amice-Vélu et Vishik (\cite{amivel}, \cite{vis}).

\begin{defi}
On appelle distribution tempérée d'ordre $r$ sur $\OF$ une forme linéaire continue sur l'espace de Banach $C^r(\OF,J',\underline{d}_{S \setminus J'})$.
\end{defi}

Notons:
\[
\big(C^r(\OF,J',\underline{d}_{S \setminus J'})^{\vee}, \|\cdot \|_{\mathcal{D}_r,J',(d_{\sigma})_{\sigma}}\big)
\] 
l'espace des distributions tempérées d'ordre $r$ sur $\OF$ muni de la topologie forte.

Soit $N\in \mathbb{Z}_{\geq 0}$. Si $\mu \in \mathcal{F}^N(\OF,J,\underline{d}_{S \setminus J})^{\vee}$ et $f \in \mathcal{F}^N(\OF,J,\underline{d}_{S \setminus J})$ on note $\int_{\OF}f(z)\mu(z)$ l'accouplement et on pose:
\[
\int_{D(a,n)} f(z)\mu(z) = \int_{\OF} \mathbf{1}_{D(a,n)}(z) f(z) \mu(z)
\] 
où, pour $a \in \OF$ et $n \in \mathbb{Z}_{\geq 0}$,  $\mathbf{1}_{D(a,n)}$ désigne la fonction caractéristique de $a+\varpi_F^n \OF$.

\begin{theo}\label{velu}  (i) Soit $\mu \in C^r(\OF,J',\underline{d}_{S \setminus J'})^{\vee}$. Il existe une constante $C_{\mu} \in \mathbb{R}_{\geq 0}$ telle que pour tout $a \in \OF$, tout $n \in \mathbb{Z}_{\geq 0}$, tout $\underline{0} \leqslant \underline{n}_{S\setminus J'} \leqslant \underline{d}_{S\setminus J'}$ et tout $\underline{m}_{J'} \in \mathbb{Z}_{\geq 0}^{|J'|}$ on ait: 
\begin{align}\label{inte1}
\Big|\int_{D(a,n)} (z-a)^{\underline{n}_{S\setminus J'}} (z-a)^{\underline{m}_{J'}} \mu(z)\Big| \leq C_{\mu} \, q^{n(r-|\underline{n}_{S\setminus J'}| - |\underline{m}_{J'}|)}. 
\end{align}

(ii) Soit $N$ un entier tel que $N\geq [r]$ et $\mu \in  \mathcal{F}^N(\OF,J,\underline{d}_{S \setminus J})^{\vee}$. Supposons qu'il existe une constante $C_{\mu} \in \mathbb{R}_{\geq 0}$ telle que pour tout $a \in \OF$, tout $n \in \mathbb{Z}_{\geq 0}$, tout $\underline{0} \leqslant \underline{n}_{S\setminus J} \leqslant \underline{d}_{S\setminus J}$ et tout $\underline{m}_{J} \in \mathbb{Z}_{\geq 0}^{|J|}$ tels que $|\underline{n}_{S\setminus J}|+|\underline{m}_{J}|\leq N$, on ait:    
\begin{align}\label{inte2}
\Big|\int_{D(a,n)} (z-a)^{\underline{n}_{S\setminus J}} (z-a)^{\underline{m}_{J}} \mu(z)\Big| \leq C_{\mu} \, q^{n(r-|\underline{n}_{S\setminus J}| - |\underline{m}_{J}|)}.
\end{align}
Alors $\mu$ se prolonge de manière unique en une distribution tempérée d'ordre $r$ sur $\OF$. 
\begin{proof}
\cite[Théorème 4.2]{marco}.
\end{proof}
\end{theo}

\begin{rema}
{\rm
La preuve du Théorème \ref{velu} utilise de manière cruciale la construction explicite d'une base de Banach pour l'espace $C^r(\OF,J',\underline{d}_{S \setminus J'})$, qui dépend de $r$ et qui consiste d'une famille dénombrable de fonctions localement polynômiales  \cite[Proposition 3.15]{marco}. Si $F = \Q$ cette base coïncide avec celle construite par  Van der Put pour l'espace des fonctions continues sur $\Z$ et généralisée par Colmez pour $r$ quelconque (\cite{vander}, \cite[Théorème I.5.14]{colmez}). Signalons que pour l'espace des fonctions continues sur $\OF$ cette base a déjà été costruite par De Shalit \cite[§2]{sha}.}
\end{rema}

\begin{rema}\label{normaequiv} 
{\rm Une conséquence facile du Théorème \ref{velu} (\cite[Corollaire 4.3]{marco}) est la remarque suivante. Si l'on définit $\|\mu \|_{r,\underline{d}_{S\setminus J}}$, pour $\mu \in C^r(\OF,J', \underline{d}_{S\setminus J'})^{\vee}$ par la formule
\[
\|\mu \|_{r,\underline{d}_{S\setminus J}} = \sup_{a \in \OF, n \in \mathbb{Z}_{\geq 0}} \sup_{\substack{ \underline{m}_J\in \mathbb{Z}_{\geq 0}^{|J|} \\ \underline{0}\leqslant  \underline{n}_{S\setminus J} \leqslant \underline{d}_{S\setminus J}}} \Big( \Big| \int_{D(a,n)} (z-a)^{\underline{n}_{S\setminus J}} (z-a)^{\underline{m}_{J}} \mu(z)\Big| q^{-n(r- |\underline{n}_{S\setminus J}| - |\underline{m}_{J}|)} \Big)
\] 
alors $\|\cdot \|_{r,\underline{d}_{S\setminus J}}$ est une norme sur $C^r(\OF,J', \underline{d}_{S\setminus J})^{\vee}$ équivalente à  $\|\cdot \|_{\mathcal{D}_r,J',(d_{\sigma})_{\sigma}}$.}
\end{rema}

\begin{rema}\label{alternativa}
{\rm
Notons $d = [F:\Q]$. En utilisant le fait que $\OF$ est un $\Z$-module libre de rang $d$ on est amené à considérer une autre notion, tout à fait naturelle, de fonction de classe $C^r$ sur $\OF$. Fixons une  $d$-uplet $\vec{r}= (r_i)_{1 \leq i \leq d}$ de nombres réels positifs ou nuls tels que $\sum r_i = r$ et une base $(e_i)_{1 \leq i \leq d}$ de $\OF$ sur $\Z$.  Notons $\theta$ l'isomorphisme de $\Z$-modules défini par:
\[
\theta\colon \Z^d \stackrel{\sim}{\longrightarrow} \OF, \quad (a_1, \ldots, a_d) \mapsto \sum_{i=1}^d a_i e_i.
\]
Si $z \in \bigotimes_{i=1}^d C^{r_i}(\Z,E)$, on définit $\|z \|$ comme l'infimum des $\sup_{j\in J} \|v_{j_1}\|_{C^{r_1}}\cdot \ldots \cdot \|v_{j_d}\|_{C^{r_d}}$ pour toutes les écritures possibles de $z$ sous la forme $\sum_{j\in J} v_{j_1}\otimes \ldots \otimes v_{j_d}$. Ceci munit $\bigotimes_{i=1}^d C^{r_i}(\Z,E)$ d'une semi-norme et on note $\widehat{\bigotimes}_{i=1}^d C^{r_i}(\Z,E)$ le séparé complété de l'espace $\bigotimes_{i=1}^d C^{r_i}(\Z,E)$ pour cette semi-norme. Notons:
\[
C^{\vec{r}}(\OF,E) = \big\{f\colon \OF \to E, \ \textstyle{f\circ \theta \in  \widehat{\bigotimes}_{i=1}^d C^{r_i}(\Z,E)}  \big\},
\]
et munissons $C^{\vec{r}}(\OF,E)$ de la topologie déduite de celle définie sur $\widehat{\bigotimes}_{i=1}^d C^{r_i}(\Z,E)$. Dans \cite[§5]{marco} on a montré que les espaces de Banach $C^r(\OF,E)$ et $C^{\vec{r}}(\OF,E)$ ne sont pas isomorphes dès que $r>0$.}
\end{rema}

\section{Représentations de $\GL(F)$}

\subsection{Généralités}

On fixe désormais une partie $J$ de $S$ jusqu'à la fin de l'article. Soit $G$ un groupe de Lie localement $F$-analytique. On note $G_0$ le groupe de Lie localement $\Q$-analytique obtenu par restriction des scalaires de $F$ à $\Q$ à partir de $G$ (\cite[§5.14]{bv}).  Si $V$ est un $E$-espace vectoriel localement convexe séparé, on peut définir, suivant \cite[§2]{sch1}, l'espace des fonctions localement $\Q$-analytiques de $G$ dans $V$ comme étant l'espace des fonctions localement analytiques de $G_0$ dans $V$. On note $\C^{\Q-an}(G,V)$ l'espace de ces fonctions muni de l'action à gauche de $G$ usuelle. 

Soit $\mathfrak{g}$ l'algèbre de Lie de $G$. On a une action $\Q$-linéaire de  $\mathfrak{g}$ sur l'espace $\C^{\Q-an}(G,V)$ définie par:
\[
(\mathfrak{x}f)(g) =  \frac{d}{dt}\Big(t \mapsto f(\exp (-t\mathfrak{x})g)\Big)\Big|_{t=0}
\]
où $\exp\colon \mathfrak{g} \dashrightarrow G$ désigne l'application exponentielle définie localement autour de $0$ \cite[§2]{sch1}. Cette action se prolonge en une action de l'algèbre de Lie $\mathfrak{g}\otimes_{\Q}E$. Comme $\mathfrak{g}$ est un $F$-espace vectoriel, alors $\mathfrak{g}\otimes_{\Q}E$ est une algèbre de Lie sur l'anneau $F\otimes_{\Q}E$. On en déduit un isomorphisme de $E$-espaces vectoriels: 
\begin{align}\label{isodilie}
\mathfrak{g}\otimes_{\Q}E \simeq  \bigoplus_{\sigma \in S} \mathfrak{g} \otimes_{F,\sigma}E.  
\end{align}

\begin{defi}[\cite{scr}, Définition 1.3.1] 
Une fonction localement $\Q$-analytique $f\colon G \to V$ est dite localement $J$-analytique si l'action de $\mathfrak{g}\otimes_{\Q}E$ sur $f$ se factorise par $\bigoplus_{\sigma \in J} \mathfrak{g} \otimes_{F,\sigma} E$.
\end{defi}

L'ensemble des fonctions localement $J$-analytiques est un sous-espace fermé de $C^{\Q-an}(G,V)$. On le munit de la topologie induite et on le note $C^{J-an}(G,V)$.

\begin{defi}[\cite{scr}, Définition 1.3.4]\label{benj}
Soit $V$ un espace vectoriel muni d'une topologie séparée localement convexe tonnelée. On dit que $V$ est une représentation localement $J$-analytique de $G$ si les deux conditions suivantes sont vérifiées: 
\begin{itemize}
\item[(i)] Le groupe $G$ agit sur $V$ par endomorphismes continus; 
\item[(ii)]  Pour tout $v \in V$, l'application de $G$ dans $V$ définie par $g \mapsto gv$ est localement $J$-analytique. 
\end{itemize}
\end{defi}

Dans la Définition \ref{benj}, l'hypothèse que $V$ soit tonnelé implique, en utilisant le Théorème de Banach-Steinhaus (\cite[Théorème 6.15]{sch}), que l'action de $G$ soit continue.

\begin{exem}
{\rm L'espace localement convexe $C^{J-an}(G,V)$ muni de l'action à gauche de $G$ usuelle est une représentation localement $J$-analytique.}
\end{exem}

\subsection{Rappels sur les induites localement analytiques de $\GL(F)$}\label{indo}

On pose $G = \GL(F)$.  On désigne par $T$ le tore déployé constitué par les matrices diagonales de $G$ et par $P$ le sous-groupe de Borel des matrices triangulaires supérieures. On désigne par $N$ le sous-groupe de $G$ des matrices unipotentes supérieures. 

Si $(\rho,P)$ est une représentation localement $J$-analytique de $P$, on note $\mathrm{Ind}_P^G(\rho)^{J-an}$ l'espace des fonctions localement $J$-analytiques de $G$ dans $V$ telles que:
\[
\forall g\in G, \forall p\in P, \quad f(pg) = \rho(p)f(g).
\] 
On munit cet espace d'une action à gauche et $E$-linéaire de $G$ par $(g f)(g') = f(g'g)$, ce qui en fait une représentation localement $J$-analytique.

Soit $\chi$ un caractère localement $\Q$-analytique de $T$. Par inflation on peut aussi le voir comme représentation localement $\Q$-analytique de $P$. Nous allons construire ici certaines  sous-représentations localement $\Q$-analytique de $\mathrm{Ind}_P^G (\chi)^{S-an}$ et puis, en utilisant l'espace des fonctions localement analytiques sur $\OF$ construit au §\ref{ana}, on en donne une description équivalente. 

Pour $t_1,t_2 \in F^{\times}$ assez proches de $1$ on peut écrire
\[
\chi([\begin{smallmatrix}t_1 & 0 \\ 0 & t_2\end{smallmatrix}]) = \prod_{\sigma\in S}  \sigma(t_1)^{d_{1,\sigma}} \sigma(t_2)^{d_{2,\sigma}},
\]
avec $d_{1,\sigma},d_{2,\sigma} \in E$. Notons $J$ le sous-ensemble des $\sigma \in S$ tels que
\[
d_{2,\sigma} - d_{1,\sigma} \notin \mathbb{Z}_{\geq 0}.
\] 
Quitte à considérer la représentation $\mathrm{Ind}_P^G (\chi)^{S-an} \otimes ((\prod_{\sigma \in S\setminus J} \sigma^{d_{1,\sigma}})\circ \mathrm{d\acute{e}t})^{-1}$, on peut supposer qu'au voisinage de $1$ on a:
\[
\chi([\begin{smallmatrix}t_1 & 0 \\ 0 & t_2\end{smallmatrix}]) = \chi_1(t_1) \chi_2(t_2) \prod_{\sigma\in S\setminus J} \sigma(t_2)^{d_{\sigma}},
\]
où $\chi_1$ et $\chi_2$ sont des caractères localement $J$-analytiques de $P$. On pose $u = [\begin{smallmatrix}0 & 0 \\ 1 & 0\end{smallmatrix}]$ et si $\sigma \in S$ on note $u_{\sigma}$ l'élément de $\mathfrak{gl}_2(F)\otimes_{\Q} E$ obtenu par l'isomorphisme \eqref{isodilie}. Si $\sigma \in S\backslash J$, on note $\mathfrak{z}_{\sigma} = (u_{\sigma})^{d_{\sigma}+1}$ et on pose:
\[
\epsilon_{\sigma}([\begin{smallmatrix}t_1 & 0 \\ 0 & t_2\end{smallmatrix}]) = \sigma(t_1 t_2^{-1}). 
\] 
D'après \cite[Proposition 1.3.11]{scr} l'élément $\mathfrak{z}_{\sigma}$ induit une application, que l'on note encore $\mathfrak{z}_{\sigma}$, de $\mathrm{Ind}_P^G (\chi)^{S-an}$ dans $\mathrm{Ind}_P^G (\chi \epsilon_{\sigma}^{d_{\sigma}})^{S-an}$ qui est surjective et dont le noyau est isomorphe à
\[
(\mathrm{Sym}^{d_{\sigma}}E^2)^{\sigma} \otimes_E \mathrm{Ind}_P^G (\chi^{\sigma})^{S\setminus\{\sigma\}-an},
\]
où
\begin{itemize} 
\item[$\bullet$] pour $\sigma \in S$ et $d_{\sigma} \in \mathbb{Z}_{\geq 0}$ on note  $(\mathrm{Sym}^{d_{\sigma}}E^2)^{\sigma}$ la représentation algébrique irréductible de $\mathrm{GL}_2 \otimes_{F,\sigma} E$ dont le plus haut poids est $\chi_{\sigma}\colon \mathrm{diag}(x_1,x_2) \mapsto \sigma(x_2)^{d_{\sigma}}$ vis-à-vis du sous-groupe des matrices triangulaires supérieurs.
\item[$\bullet$] On désigne par $\chi^{\sigma}$ le caractère
\[
\chi_1 \otimes \chi_2 \prod_{\tau \in S\setminus (J \coprod \{\sigma\})} \tau^{d_{\tau}}.
\]
\end{itemize}
On en déduit immédiatement pour toute partie $S'$ de $S\backslash J$ l'isomorphisme suivant:
\[
\bigcap_{\sigma \in S'} \mathrm{ker}\, \mathfrak{z}_{\sigma} \stackrel{\sim}{\longrightarrow}
 \Big(  \bigotimes_{\sigma \in S'} (\Sym^{d_{\sigma}} E^2)^{\sigma}  \Big) \otimes_E \Big(\Ind_{P}^{G} \chi_1 \otimes \chi_2 \prod_{(S\setminus J)\setminus S'} \sigma^{d_{\sigma}} \Big)^{S\setminus S'-{an}}.
\]
Notons $m_{\sigma} = d_{\sigma}+1$. D'après la preuve de \cite[Proposition 1.3.11]{scr} on a le diagramme commutatif suivant:
\[
\xymatrix@R+12pt@C+42pt{
\mathrm{Ind}_P^G (\chi)^{S-an} \ar[r]^{\mathfrak{z}_{\sigma}} \ar[d]_{}
&  \mathrm{Ind}_P^G (\chi \epsilon_{\sigma}^{d_{\sigma}})^{S-an}  \ar[d]_{} \\
 (\mathcal{F}(\OF,S))^2 \ar[r]^{\phantom{}\Big(-\frac{\partial^{m_{\sigma}}}{\partial z_{\sigma}^{m_{\sigma}}}, -\frac{\partial^{m_{\sigma}}}{\partial z_{\sigma}^{m_{\sigma}}} \Big)} & (\mathcal{F}(\OF,S))^2}   
\]
où
\begin{itemize}
\item[$\bullet$]  $\mathcal{F}(\OF,S)$ désigne l'espace $\mathcal{F}(U,J,\underline{d}_{S\setminus J})$ pour $U = \OF$ et $J=S$ (et donc $S\backslash J = \emptyset$); 
\item[$\bullet$] l'application verticale de gauche (resp. de droite) est un isomorphisme topologique explicitement donné par:
\begin{align*}
f  \longmapsto \Big((z\mapsto f ([\begin{smallmatrix}0 & 1 \\ -1 & \varpi_F z\end{smallmatrix}])), (z\mapsto f ([\begin{smallmatrix}1 & 0 \\ z & -1\end{smallmatrix}]) )\Big). 
\end{align*}
\end{itemize}
On en déduit un isomorphisme topologique:
\begin{align}\label{traduiso}
\Big(  \bigotimes_{\sigma \in S'} (\Sym^{d_{\sigma}} E^2)^{\sigma}  \Big) \otimes_E \Big(\Ind_{P}^{G} \chi_1 \otimes \chi_2 \prod_{(S\setminus J)\setminus S'} \sigma^{d_{\sigma}} \Big)^{S\setminus S'-{an}} \simeq (\mathcal{F}(\OF,S\backslash S',\underline{d}_{S'}))^2 
\end{align}
Posons:
\[
I(\chi,S\backslash S',\underline{d}_{S'}) = \Big(  \bigotimes_{\sigma \in S'} (\Sym^{d_{\sigma}} E^2)^{\sigma}  \Big) \otimes_E \Big(\Ind_{P}^{G} \chi_1 \otimes \chi_2 \prod_{(S\setminus J)\setminus S'} \sigma^{d_{\sigma}} \Big)^{S\setminus S'-{an}}
\]
et notons $V$ le $E$-espace vectoriel des fonctions $f\colon F \to E$ qui vérifient les deux conditions suivantes:  
\begin{itemize}
\item[(i)] $f|_{\OF}$ est dans $\mathcal{F}(\OF,S\backslash S',\underline{d}_{S'})$;
\item[(ii)] $\chi_2 \chi_1^{-1}(z) z^{\underline{d}_{S\setminus J}}   f(1/z)|_{\OF-\{0\}}$ se prolonge sur $\OF$ en une fonction dans $\mathcal{F}(\OF,S\backslash S',\underline{d}_{S'})$.
\end{itemize} 
L'application:
\begin{equation}\label{traduiso2}
\begin{aligned}
V &\longrightarrow \quad \mathcal{F}(\OF,S\backslash S',\underline{d}_{S'})\oplus \mathcal{F}(\OF,S\backslash S',\underline{d}_{S'})  \\ 
f &\longmapsto \Big( \big( z \mapsto f(\varpi_F z)\big), \big(z \mapsto \chi_2 \chi_1^{-1}(z) z^{\underline{d}_{S\setminus J}}   f(1/z) \big) \Big) 
\end{aligned}
\end{equation}
est un isomorphisme de $E$-espaces vectoriels. On munit $V$ de la topologie localement convexe déduite de cette application. Par les isomorphismes \eqref{traduiso} et \eqref{traduiso2} et d'après  l'égalité
\[
\begin{bmatrix} {0} & {1} \cr {-1} & {z} \end{bmatrix} \begin{bmatrix} {a} & {b} \cr {c} & {d} \end{bmatrix} = \begin{bmatrix} {\frac{ad-bc}{-cz+a}} & {-c} \cr {0} & {-cz+a} \end{bmatrix} \begin{bmatrix} {0} & {1} \cr {-1} & {\frac{dz-b}{-cz+a}} \end{bmatrix}
\]
on déduit que l'action de $G$ sur $I(\chi,S\backslash S',\underline{d}_{S'})$ se traduit sur $V$, pour tout $g = [\begin{smallmatrix}a & b \\ c & d\end{smallmatrix}] \in G$ et tout $f\in V$, par la formule
\begin{align}\label{azioneanalitica}
\left(\begin{bmatrix} {a} & {b} \cr {c} & {d} \end{bmatrix}  f \right)(z) = \chi_1(\mathrm{d\acute{e}t}(g)) \chi_2 \chi_1^{-1}(-cz+a) (-cz+a)^{\underline{d}_{S\setminus J}} f\left(\frac{d z-b}{-cz+a}\right)
\end{align}
pour tout $z\in F$, $z\neq \frac{a}{c}$, et que l'on peut prolonger $g f$ par continuité en $z = \frac{a}{c}$ (si $c \neq 0$) en une fonction appartenant à $V$.


\subsection{Une $\GL(F)$-représentation de Banach}\label{repban}
Soit $J\subseteq S$, $\chi_1, \chi_2\colon F^{\times} \to E^{\times}$ deux caractères localement $J$-analytiques et $\underline{d}_{S\setminus J}$ un $|S\backslash J|$-uplet d'entiers positifs ou nuls. Notons $r = - \mathrm{val}_{\Q}(\chi_1(p))$ et supposons $r \geq 0$. Posons:
\begin{align*}
J' &= J \coprod \{\sigma \in S\backslash J,\, d_{\sigma}+1 > r\},\quad \chi_1' = \chi_1,\quad
\chi_2' = \chi_2 \prod_{\sigma\in J'\setminus J}\sigma^{d_{\sigma}}.
\end{align*}

Nous allons ici construire un  $G$-Banach (Définition \ref{banachunit}) en utilisant les espaces qui ont été définis au §\ref{chiusi}.

Notons $B(\chi',J',\underline{d}_{S\setminus J'})$ le $E$-espace vectoriel des fonctions $f\colon F \to E$ qui vérifient les deux conditions suivantes: 
\begin{itemize}
\item[(i)] $f|_{\OF}$ est une fonction dans $C^{r}(\OF,J',\underline{d}_{S\setminus J'})$;
\item[(ii)] $\chi_2' {\chi_1'}^{-1}(z) z^{\underline{d}_{S\setminus J'}} f(1/z)|_{\OF-\{0\}}$ se prolonge sur $\OF$ en une fonction dans $C^{r}(\OF,J',\underline{d}_{S\setminus J'})$.
\end{itemize} 

L'application:
\begin{equation}\label{isocr}
\begin{aligned}
B(\chi',J',\underline{d}_{S\setminus J'}) &\longrightarrow \quad C^r(\OF,J',\underline{d}_{S\setminus J'})\oplus C^r(\OF,J',\underline{d}_{S\setminus J'})  \\ 
f &\longmapsto \Big( \big( z \mapsto f(\varpi_F z)\big), \big(z \mapsto \chi_2' {\chi_1'}^{-1}(z) z^{\underline{d}_{S\setminus J'}}   f(1/z) \big) \Big) 
\end{aligned}
\end{equation}
est un isomorphisme de $E$-espaces vectoriels. On munit $B(\chi',J',\underline{d}_{S\setminus J'})$ de la topologie localement convexe déduite de cette application, ce qui en fait un espace de Banach sur $E$. Plus précisement, si on désigne par $(f_1,f_2)$ l'élément de $(C^{r}(\OF,J',\underline{d}_{S\setminus J'}))^2$ qui correspond à $f \in B(\chi',J',\underline{d}_{S\setminus J'})$ via l'isomorphisme \eqref{isocr}, on a:
\begin{align}
\|f \|_B = \sup \big(\|f_1\|_{C^r}, \|f_2\|_{C^r}\big).
\end{align}
Pour $f\in B(\chi',J',\underline{d}_{S\setminus J'})$ et $g = [\begin{smallmatrix}a & b \\ c & d\end{smallmatrix}] \in G$ considérons la fonction définie par:
\begin{align}\label{acb}
\left(\begin{bmatrix} {a} & {b} \cr {c} & {d} \end{bmatrix}  f \right)(z) = \chi_1(\mathrm{d\acute{e}t}(g)) \chi_2' {\chi_1'}^{-1}(-cz+a) (-cz+a)^{\underline{d}_{S\setminus J'}} f\left(\frac{d z-b}{-cz+a}\right)
\end{align}
pour tout $z \neq \frac{a}{c}$ (si $c \neq 0$). Le résultat suivant montre que $g f$ se prolonge par continuité en $z = \frac{a}{c}$ en une fonction appartenant à $B(\chi',J',\underline{d}_{S\setminus J'})$ et que, muni de l'action de $G$ définie par:
\[
(g,f) \mapsto g f,
\]
l'espace $B(\chi',J',\underline{d}_{S\setminus J'})$ devient un $G$-Banach.

\begin{lemm}\label{autocont}
L'action à gauche de $G$ sur l'espace $B(\chi',J',\underline{d}_{S\setminus J'})$ décrite par la formule \eqref{acb} est bien définie et se fait par automorphismes continus.
\begin{proof}

Soit $f = (f_1,f_2) \in B(\chi',J',\underline{d}_{S\setminus J'})$. En utilisant l'isomorphisme \eqref{isocr} il est facile de voir que pour tout $g = [\begin{smallmatrix}a & b \\ c & d\end{smallmatrix}] \in G$ on a:
\begin{align*}
(g f)_1 (z) = \chi_1'( \mathrm{det}(g)) \chi_2' {\chi_1'}^{-1}(-c \varpi_F z+a) (-c\varpi_F z+a)^{\underline{d}_{S\setminus J'}} f_1 \Big( \frac{d z-\frac{b}{\varpi_F}}{-c\varpi_F z+a}  \Big)
\end{align*}
si $\frac{d\varpi_F z-b}{-c\varpi_F z+a} \in \varpi_F \OF$ et 
\begin{align*}
(g f)_1 (z) = \chi_1'(\mathrm{det}(g)) \chi_2' {\chi_1'}^{-1}(d\varpi_F z-b)  (d \varpi_F z-b)^{\underline{d}_{S\setminus J'}} f_2 \Big( \frac{-c\varpi_F z+a}{d\varpi_F z-b}  \Big)
\end{align*}
si $\frac{d \varpi_F z-b}{-c \varpi_F z+a} \in F\backslash \varpi_F \OF$;
\begin{align*}
(g f)_2 (z) = \chi_1'(\mathrm{det}(g)) \chi_2' {\chi_1'}^{-1}(-c +az) (-c+az)^{\underline{d}_{S\setminus J'}} f_1 \Big( \frac{-b\frac{z}{\varpi_F}+ \frac{d}{\varpi_F}}{a z-c}  \Big)
\end{align*}
si $\frac{-bz+d}{a z-c} \in \varpi_F \OF$ et 
\begin{align*}
(g f)_2 (z) = \chi_1'( \mathrm{det}(g)) \chi_2' {\chi_1'}^{-1}(-bz +d) (-bz+d)^{\underline{d}_{S\setminus J'}} f_2 \Big( \frac{az-c}{-bz+d}  \Big)
\end{align*}
si $\frac{-bz+d}{a z-c} \in F\backslash \varpi_F \OF$.

Il suffit donc de montrer que l'application
\begin{equation}\label{contin}
\begin{aligned}
C^{r}(\OF,J',\underline{d}_{S\setminus J'}) \oplus C^{r}(\OF,J',\underline{d}_{S\setminus J'}) &\longrightarrow C^r(\OF,J',\underline{d}_{S\setminus J'})\oplus C^r(\OF,J',\underline{d}_{S\setminus J'})  \\ 
(f_1, f_2) \quad \quad \quad\quad \quad\quad\ &\longmapsto \quad\quad\quad\quad\  ((g f)_1,(g f)_2)
\end{aligned}
\end{equation}
est bien définie et continue. Rappelons que par la décomposition de Bruhat on a:
\begin{align}\label{bruhat}
G = P \cup PwN.
\end{align}
On est alors réduit à montrer la stabilité et la continuité de l'application \eqref{contin} pour les matrices $g$ de la forme $[\begin{smallmatrix}\lambda & 0 \\ 0 & \lambda \end{smallmatrix}]$, $[\begin{smallmatrix}0 & \varpi_F \\ 1 & 0\end{smallmatrix}]$,  $[\begin{smallmatrix}1 & 0 \\ 0 & \lambda\end{smallmatrix}]$ et $[\begin{smallmatrix}1 & \lambda \\ 0 & 1\end{smallmatrix}]$ (avec $\lambda \in F^{\times}$). Or ceci est une conséquence des formules ci-dessus, de la Proposition \ref{composizione} et du fait que l'espace $C^{r}(\OF,J',\underline{d}_{S\setminus J'})$ est une $E$-algèbre de Banach   (\cite[Lemme 2.9]{marco}).

\end{proof}
\end{lemm}

\begin{rema}
{\rm Le Lemme \ref{autocont} et le Théorème de Banach-Steinhaus \cite[Théorème 6.15]{sch} impliquent que l'espace $B(\chi',J',\underline{d}_{S\setminus J'})$ est un $G$-Banach.}
\end{rema}

Soit $k \in \mathbb{Z}_{>0}$. Notons $S_k \subset \OF^{\times}$ un système de représentants des classes de $(\OF/\varpi_F^{k}\OF)^{\times}$. Notons $l$ le plus petit entier positif tel que ${\chi_1'}|_{D(a_i,l)}$ (resp. \ ${\chi_2'}|_{D(a_i,l)}$) 
est une fonction $J'$-analytique sur l'ouvert $D(a_i,l)$ pour tout $a_i \in S_{l}$.

Faisons l'hypothèse supplémentaire suivante:
\begin{align}\label{iposup}
\mathrm{val}_{\Q}(\chi_1(p)) + \mathrm{val}_{\Q}(\chi_2(p)) + |\underline{d}_{S\setminus J}| = 0,
\end{align}
et notons que \eqref{iposup} est équivalente à
\begin{align}\label{iposup33}
\mathrm{val}_{\Q}(\chi_1'(p)) + \mathrm{val}_{\Q}(\chi_2'(p)) + |\underline{d}_{S\setminus J'}| = 0.
\end{align}
\begin{lemm}\label{funzionicr}
Les fonctions de $F$ dans $E$: 
\begin{align*}
z &\mapsto z^{\underline{n}_{S\setminus J'}} z^{\underline{m}_{J'}} \\
z &\mapsto \left\{ \begin{array}{ll}
\chi_2'{\chi_1'}^{-1}(z-a) (z-a)^{\underline{d}_{S\setminus J'}- \underline{n}_{S\setminus J'}} (z-a)^{-\underline{m}_{J'}}   & \mbox{si} \ z \neq a \\
0 & \mbox{si} \ z=a.
\end{array} \right.
\end{align*} 
pour tout $a \in F$, tout $\underline{m}_{J'} \in \mathbb{Z}_{\geq 0}^{|J'|}$ et tout $\underline{0}\leqslant \underline{n}_{S\setminus J'} \leqslant \underline{d}_{S\setminus J'}$ tels que $r - \big(|\underline{n}_{S\setminus J'}| + |\underline{m}_{J'}|\big) > 0$ sont dans $B(\chi',J',\underline{d}_{S\setminus J'})$. 
\begin{proof}
Le même raisonnement que dans \cite[Lemme 4.2.2]{bb} s'applique. Il suffit de montrer que la fonction de $\OF$ dans $E$ définie par:
\[
z \mapsto \left\{ \begin{array}{ll}
\chi_2'{\chi_1'}^{-1}(z) z^{\underline{d}_{S\setminus J'}- \underline{n}_{S\setminus J'}} z^{-\underline{m}_{J'}}   & \mbox{si} \ z \neq 0 \\
0 & \mbox{si} \ z=0
\end{array} \right.
\]
est dans $C^{r}(\OF,J',\underline{d}_{S\setminus J'})$. Soit $f_0$ la fonction nulle sur $\OF$ et, pour $n \in \mathbb{Z}_{> 0}$ posons:
\[
f_n(z) = \mathbf{1}_{\OF \setminus D(0,n)}(z)\chi_2'{\chi_1'}^{-1}(z) z^{\underline{d}_{S\setminus J'}- \underline{n}_{S\setminus J'}} z^{-\underline{m}_{J'}}.
\] 
La fonction $f_n$ est bien dans $C^{r}(\OF,J',\underline{d}_{S\setminus J'})$ puisqu'elle est en particulier dans $\mathcal{F}(\OF, J',\underline{d}_{S\setminus J'})$. Par \cite[Lemme 9.9]{sch} il suffit de montrer que $f_{n+1}-f_n$ tend vers $0$ dans l'espace de Banach dual de l'espace de Banach des distributions tempérées d'ordre $r$ sur $\OF$, i.e.:
\[
\sup_{\mu \in C^r(\OF,J',\underline{d}_{S \setminus J'})^{\vee}} \frac{\Big|\int_{\OF}\Big(f_{n+1}(z)- f_{n}(z)\Big)\mu(z)\Big|}{\|\mu \|_{r,\underline{d}_{S \setminus J}}} \to 0 \ \mbox{quand} \ n \to +\infty.
\]
Notons que l'on a:
\begin{equation}\label{rapmod}
\begin{aligned}
f_{n+1}(z)-f_{n}(z) &= \mathbf{1}_{D(0,n)\setminus D(0,n+1)}(z)\chi_2'{\chi_1'}^{-1}(z) z^{\underline{d}_{S\setminus J'}- \underline{n}_{S\setminus J'}}z^{-\underline{m}_{J'}} \\
&= \sum_{a_i \in S_{l}} \mathbf{1}_{D(a_i \varpi_F^n,n+l)}(z) \chi_2'{\chi_1'}^{-1}(z) z^{\underline{d}_{S\setminus J'}- \underline{n}_{S\setminus J'}}z^{-\underline{m}_{J'}}.
\end{aligned}
\end{equation}
Comme $\chi_1'$ et $\chi_2'$ sont des caractères $J'$-analytiques sur $D(a_i,l)$ pour tout $a_i \in S_{l}$ on a pour tout $n \geq 0$:
\begin{align*}
 \mathbf{1}_{D(a_i \varpi_F^n,n+l)}(z) \chi_2'{\chi_1'}^{-1}(z) &= \chi_2'{\chi_1'}^{-1}(\varpi_F^n)   \mathbf{1}_{D(a_i,l)}\Big(\frac{z}{\varpi_F^n}\Big) \chi_2'{\chi_1'}^{-1}\Big( \frac{z}{\varpi_F^n} \Big) \\
&=\chi_2'{\chi_1'}^{-1}(\varpi_F^n)  \mathbf{1}_{D(a_i,l)}\Big(\frac{z}{\varpi_F^n} \Big) \sum_{\underline{h}_{J'} \geqslant \underline{0}} b_{\underline{h}_{J'}}(a_i)\Big(\frac{z}{\varpi_F^n}-a_i \Big)^{\underline{h}_{J'}} \\
&= \chi_2'{\chi_1'}^{-1}(\varpi_F^n)  \sum_{\underline{h}_{J'} \geqslant \underline{0}}  \mathbf{1}_{D(a_i \varpi_F^n,n+l)}(z)  b_{\underline{h}_{J'}}(a_i)\Big(\frac{z-a_i\varpi_F^n}{\varpi_F^n}  \Big)^{\underline{h}_{J'}}. 
\end{align*}
Notons $C_1 = \sup_{a_i \in S_{l}} \sup_{\underline{h}_{J'}} |b_{\underline{h}_{J'}}(a_i)|$ et remarquons que la condition \eqref{iposup33} implique l'égalité $\big|\chi_2' {\chi_1'}^{-1}(\varpi_F^n)  \big| = q^{-n(2r-|\underline{d}_{S\setminus J'}|)}$.

En écrivant $z^{-\underline{m}_{J'}} = (z-a_i\varpi_F^n+a_i\varpi_F^n)^{-\underline{m}_{J'}}$ et en développant on obtient pour tout $a_i \in S_{l}$:
\[
\mathbf{1}_{D(a_i \varpi_F^n,n+l)}(z)  z^{-\underline{m}_{J'}} = \mathbf{1}_{D(a_i \varpi_F^n,n+l)}(z) (a_i\varpi_F^n)^{-\underline{m}_{J'}} \sum_{\underline{t}_{J'} \geqslant \underline{0}} \lambda_{\underline{t}_{J'}} a_i^{-\underline{t}_{J'}} \Big(\frac{z-a_i\varpi_F^n}{\varpi_F^n} \Big)^{\underline{t}_{J'}} 
\]
où les $\lambda_{\underline{t}_{J'}}$ sont des éléments de $\Oe$. De manière analogue, on obtient pour tout $a_i \in S_{l}$: 
\begin{align*}
 & \ \mathbf{1}_{D(a_i \varpi_F^n,n+l)}(z) z^{\underline{d}_{S\setminus J'}-\underline{n}_{S\setminus J'}} \\
=& \   \mathbf{1}_{D(a_i \varpi_F^n,n+l)}(z)\sum_{\underline{0}\leqslant \underline{k}_{S\setminus J'}\leqslant \underline{d}_{S\setminus J'}-\underline{n}_{S\setminus J'}}\mu_{\underline{k}_{S\setminus J'}} (a_i \varpi_F^n)^{\underline{k}_{S\setminus J'}} (z-a_i \varpi_F^n)^{\underline{d}_{S\setminus J'}-\underline{n}_{S\setminus J'}-\underline{k}_{S\setminus J'}},
\end{align*}
où les $\mu_{\underline{k}_{S\setminus J'}}$ sont des entiers.  

Notons $f_{\underline{\alpha}_{S\setminus J'},\underline{\beta}_{J'}}$, pour tout $\underline{0}\leqslant \underline{\alpha}_{S\setminus J'}\leqslant \underline{d}_{S\setminus J'}$ et tout $\underline{\beta}_{J'} \in \mathbb{Z}_{\geq 0}^{|J'|}$ la fonction de $\OF-\{0\}$ dans $E$ définie par:
\[
z \mapsto z^{\underline{d}_{S\setminus J'}-\underline{\alpha}_{S\setminus J'}}z^{-\underline{\beta}_{J'}}.
\]
Par \eqref{rapmod} on a:
\[
\big|\mu \big(f_{n+1}(z)-f_n(z)\big)\big| = \sup_{a_i \in S_{l}}\big|\mu(\mathbf{1}_{D(a_i \varpi_F^n,n+l)}(z) \chi_2'{\chi_1'}^{-1}(z) f_{\underline{n}_{S\setminus J'},\underline{m}_{J'}}(z)) \big|,
\]
et, en utilisant les égalités précédentes on déduit pour tout $a_i \in S_{l}$:
\begin{align*}
&\big|\mu(\mathbf{1}_{D(a_i \varpi_F^n,n+l)}(z) \chi_2'{\chi_1'}^{-1}(z) z^{\underline{d}_{S\setminus J'}- \underline{n}_{S\setminus J'}}z^{-\underline{m}_{J'}}) \big|\\
\leq & C_1 q^{-n(2r - |\underline{d}_{S\setminus J'}|-|\underline{m}_{J'}|)}  \sup_{\substack{
 \underline{l}_{J'} \\
 \underline{k}_{S\setminus J'} }}  q^{-n(|\underline{k}_{S\setminus J'}|-|\underline{l}_{J'}|)} \big|\mu\big(\mathbf{1}_{D(a_i \varpi_F^n,n+l)}(z)  f_{\underline{n}_{S\setminus J'}+\underline{k}_{S\setminus J'},\underline{l}_{J'}}(z-a_i\varpi_F^n) \big) \big|,
\end{align*} 
où $\underline{l}_{J'}$ varie dans $\mathbb{Z}_{\geq 0}^{|J'|}$ et $\underline{0}\leqslant \underline{k}_{S\setminus J'}\leqslant \underline{d}_{S\setminus J'}$. D'après la Remarque \ref{normaequiv} on a:
\[
\big|\mu\big(\mathbf{1}_{D(a_i \varpi_F^n,n+l)}(z)  f_{\underline{n}_{S\setminus J'}+\underline{k}_{S\setminus J'},\underline{l}_{J'}}(z-a_i\varpi_F^n) \big) \big| \leq \|\mu \|_{r,\underline{d}_{S \setminus J}} \sup_{\substack{
 \underline{l}_{J'} \\
 \underline{k}_{S\setminus J'} }} q^{(n+l)(r +|\underline{k}_{S\setminus J'}|-|\underline{l}_{J'}|-|\underline{d}_{S\setminus J'}|+|\underline{n}_{S\setminus J'}|)},
\]  
d'où
\[
\big|\mu \big(f_{n+1}(z)-f_n(z)\big)\big| \leq C_1 \|\mu \|_{r,\underline{d}_{S \setminus J}} q^{-n(r-|\underline{m}_{J'}|-|\underline{n}_{S\setminus J'}|)} \sup_{\substack{
 \underline{l}_{J'} \\
 \underline{k}_{S\setminus J'} }} q^{l(r +|\underline{k}_{S\setminus J'}|-|\underline{l}_{J'}|-|\underline{d}_{S\setminus J'}|+|\underline{n}_{S\setminus J'}|)}. 
\]
On en déduit le résultat car $r > |\underline{m}_{J'}|+|\underline{n}_{S\setminus J'}|$.
\end{proof}
\end{lemm}

D'après le Lemme \ref{funzionicr} on sait que pour tout $a \in F$, tout $\underline{m}_{J'} \in \mathbb{Z}_{\geq 0}^{|J'|}$ et tout $\underline{0}\leqslant \underline{n}_{S\setminus J'} \leqslant \underline{d}_{S\setminus J'}$ tels que $r - |\underline{n}_{S\setminus J'}|- |\underline{m}_{J'}| > 0$, les fonctions $z^{\underline{n}_{S\setminus J'}} z^{\underline{m}_{J'}}$ et $\chi_2'{\chi_1'}^{-1}(z-a) (z-a)^{\underline{d}_{S\setminus J'}- \underline{n}_{S\setminus J'}} (z-a)^{-\underline{m}_{J'}}$ sont dans $B(\chi',J',\underline{d}_{S\setminus J'})$. Notons $L(\chi',J',\underline{d}_{S\setminus J'})$ l'adhérence dans $B(\chi',J',\underline{d}_{S\setminus J'})$ du sous-$E$-espace vectoriel engendré par ces fonctions.

\begin{lemm}\label{stabilità}
Le sous-espace $L(\chi',J',\underline{d}_{S\setminus J'})$ est stable par $G$ dans $B(\chi',J',\underline{d}_{S\setminus J'})$.
\begin{proof}
Il s'agit d'un calcul facile et est laissé au lecteur. 
\end{proof}
\end{lemm}

Posons:
\[
\Pi(\chi',J',\underline{d}_{S\setminus J'}) \ugu B(\chi',J',\underline{d}_{S\setminus J'})/L(\chi',J',\underline{d}_{S\setminus J'}). 
\]
C'est un espace de Banach sur $E$ et, d'après les Lemmes \ref{autocont} et \ref{stabilità},  il est  muni d'une action de $G$ par automorphismes continus.

\section{Réseaux}

\subsection{Deux conditions nécessaires de non nullité}\label{lattiprov}

Soit $J\subseteq S$, $\chi_1, \chi_2\colon F^{\times} \to E^{\times}$ deux caractères localement $J$-analytiques et $\underline{d}_{S\setminus J}$ un $|S\backslash J|$-uplet d'entiers positifs ou nuls. Notons $r = - \mathrm{val}_{\Q}(\chi_1(p))$ et considérons la représentation localement $\Q$-analytique:
\[
I(\chi,J,\underline{d}_{S\setminus J}) = \Big(  \bigotimes_{\sigma \in S \setminus J} (\Sym^{d_{\sigma}} E^2)^{\sigma}  \Big) \otimes_E \Big(\Ind_{P}^{G} \chi_1 \otimes \chi_2 \Big)^{J-{an}}
\]
qui a été construite au §\ref{indo}. Soit $I(\chi,J,\underline{d}_{S\setminus J})(F)$ le sous-espace fermé de $I(\chi,J,\underline{d}_{S\setminus J})$ des fonctions $f$ qui sont à support compact. Il est stable par $P$ et il engendre $I(\chi,J,\underline{d}_{S\setminus J})$ sous $G$. En particulier, cet espace contient l'espace $\mathcal{O}(\OF,J,\underline{d}_{S\setminus J})$ et l'on peut voir facilement que:
\[
I(\chi,J,\underline{d}_{S\setminus J}) = \sum_{g\in G} g \mathcal{O}(\OF,J,\underline{d}_{S\setminus J}).
\] 
D'après la preuve de \cite[Proposition 1.21]{eme}, le complété unitaire universel de $I(\chi,J,\underline{d}_{S\setminus J})$ est le complété par rapport au sous-$\Oe[G]$-réseau engendré par les vecteurs $\mathbf{1}_{\OF}(z)z^{\underline{n}_{S\setminus J}}z^{\underline{m}_{J}}$ avec $\underline{0}\leqslant\underline{n}_{S\setminus J} \leqslant \underline{d}_{S\setminus J}$ et $\underline{m}_{J} \in \mathbb{Z}_{\geq 0}^{|J|}$. En utilisant $G = P K$ et le fait que $K$ est compact on voit qu'il suffit de compléter par rapport au sous-$\Oe[P]$-réseau engendré par les vecteurs $\mathbf{1}_{\OF}(z)z^{\underline{n}_{S\setminus J}}z^{\underline{m}_{J}}$ et $\mathbf{1}_{F-\OF}(z)\chi_2\chi_1^{-1}(z) z^{\underline{d}_{S\setminus J}-\underline{n}_{S\setminus J}}z^{-\underline{m}_{J}}$ avec $\underline{0}\leqslant\underline{n}_{S\setminus J} \leqslant \underline{d}_{S\setminus J}$ et $\underline{m}_{J} \in \mathbb{Z}_{\geq 0}^{|J|}$. Notons $\Lambda$ ce réseau et $I(\chi,J,\underline{d}_{S\setminus J} )^{\bigwedge}$ le complété de $I(\chi,J,\underline{d}_{S\setminus J})$ par rapport à $\Lambda$. C'est en particulier un $G$-Banach unitaire.

\begin{prop} \label{nullita}
Le deux conditions suivantes sont nécessaires pour que   $I(\chi,J,\underline{d}_{S\setminus J} )^{\bigwedge}$ soit non nul: 
\begin{itemize}
\item[(i)] Le caractère central de $I(\chi,J,\underline{d}_{S\setminus J} )$ est intègre; 
\item[(ii)] On a l'inégalité $\mathrm{val}_{\Q}(\chi_2(p)) + |\underline{d}_{S\setminus J}| \geq 0$.  
\end{itemize}
\begin{proof}
Supposons que $(I(\chi,J,\underline{d}_{S\setminus J} )^{\bigwedge},\|\cdot\|)$ soit non nul. En particulier l'application canonique
\[
\iota\colon I(\chi,J,\underline{d}_{S\setminus J} ) \to I(\chi,J,\underline{d}_{S\setminus J} )^{\bigwedge}
\]
est non nulle. Soit $f\in I(\chi,J,\underline{d}_{S\setminus J} )$ tel que $\iota(f) \neq 0$. Alors, comme $\iota$ est $G$-équivariante et $I(\chi,J,\underline{d}_{S\setminus J} )^{\bigwedge}$ est un $G$-Banach unitaire on a:
\[
\Big|\chi_1(p)\chi_2(p) p^{|\underline{d}_{S\setminus J}|}\Big| \|\iota(f) \| = \|\iota(f) \|, 
\]
d'où (i).

Montrons maintenant que si $\mathrm{val}_{\Q}(\chi_2(p)) + |\underline{d}_{S\setminus J}| < 0$ alors $I(\chi,J,\underline{d}_{S\setminus J} )^{\bigwedge}$ est nul. Cela est équivalent à prouver que pour tout $\underline{0}\leqslant\underline{n}_{S\setminus J} \leqslant \underline{d}_{S\setminus J}$ et tout  $\underline{m}_{J} \in \mathbb{Z}_{\geq 0}^{|J|}$ on a:
\begin{align}\label{nullita1}
\forall \lambda \in E, \forall n\geq 0, \quad \lambda \mathbf{1}_{D(0,n)}(z) z^{\underline{n}_{S\setminus J}}z^{\underline{m}_{J}} \in \Lambda. 
\end{align}  
La démonstration se fait par récurrence sur $|\underline{n}_{S\setminus J}|+|\underline{m}_{J}|$. 

Supposons $|\underline{n}_{S\setminus J}|+|\underline{m}_{J}| = 0$. Soit $\lambda \in E$ et $n \in \mathbb{Z}_{\geq 0}$. Notons $m$ le plus petit entier positif tel que $val_F(\chi_2(\uni^m)\uni^{m \underline{d}_{S\setminus J}})< val_F(\lambda)$ et fixons $R\subset \OF$ un système de représentants des classes de $\OF/\uni^m \OF$. D'après la formule \eqref{azioneanalitica} et comme $\Lambda$ est stable sous l'action de $P$ on a:
\[
\forall a_i \in R, \quad [\begin{smallmatrix}\uni^m & \uni^n a_i \\ 0 & 1\end{smallmatrix}] \mathbf{1}_{D(0,n)} =  \chi_2(\uni^m)\uni^{m \underline{d}_{S\setminus J}} \mathbf{1}_{D(\uni^n a_i, n+m)} \in \Lambda.
\]  
On en déduit:
\[
\sum_{a_i \in R} \chi_2(\uni^m)\uni^{m \underline{d}_{S\setminus J}} \mathbf{1}_{D(\uni^n a_i, n+m)} = \chi_2(\uni^m)\uni^{m \underline{d}_{S\setminus J}} \mathbf{1}_{D(0,n)} \in \Lambda,
\]
d'où $\lambda \mathbf{1}_{D(0,n)} \in \Lambda$.

Supposons que \eqref{nullita1} soit vrai pour tout   $\underline{0}\leqslant\underline{n}_{S\setminus J} \leqslant \underline{d}_{S\setminus J}$ et tout $\underline{m}_{J} \in \mathbb{Z}_{\geq 0}^{|J|}$ tels que $|\underline{n}_{S\setminus J}|+|\underline{m}_{J}|  \leq l$ où $l$ est un entier positif. Soit $\underline{i} \in \mathbb{Z}_{\geq 0}^{|S|}$ tel que:
\[
|\underline{i}| = l+1 \quad \mbox{et} \  \quad i_{\sigma}\leq d_{\sigma}, \quad \forall \sigma\in S\backslash J.  
\] 
D'après la formule \eqref{azioneanalitica} et comme $\Lambda$ est stable sous l'action de $P$ on a:
\[
\forall a_i \in R, \quad [\begin{smallmatrix}\uni^m & \uni^n a_i \\ 0 & 1\end{smallmatrix}] z^{\underline{i}}\mathbf{1}_{D(0,n)} =  \chi_2(\uni^m)\uni^{m \underline{d}_{S\setminus J}}  \Big(\frac{z-a_i \uni^n}{\uni^m} \Big)^{\underline{i}}  \mathbf{1}_{D(\uni^n a_i, n+m)} \in \Lambda,
\]  
où les $\mu_{\underline{k}}$ sont des entiers. On en déduit, en développant $\big(\frac{z-a_i \uni^n}{\uni^m} \big)^{\underline{i}}$ et en utilisant l'hypothèse de récurrence:
\begin{align}\label{duecondi}
\forall a_i \in R,\quad \chi_2(\uni^m)\uni^{m \underline{d}_{S\setminus J}} \Big(\frac{z}{\uni^m} \Big)^{\underline{i}} \mathbf{1}_{D(\uni^n a_i, n+m)} \in \Lambda.
\end{align}
En particulier, par \eqref{duecondi} on a:
\[
\sum_{a_i \in R} \chi_2(\uni^m)\uni^{m \underline{d}_{S\setminus J}} \uni^{-m\underline{i}} z^{\underline{i}} \mathbf{1}_{D(\uni^n a_i, n+m)} = \chi_2(\uni^m)\uni^{m \underline{d}_{S\setminus J}} \uni^{-m\underline{i}} z^{\underline{i}} \mathbf{1}_{D(0, n)} \in \Lambda,
\]
d'où $\lambda z^{\underline{i}} \mathbf{1}_{D(0, n)} \in \Lambda$, ce qui permet de conclure.

\end{proof}
\end{prop}

\begin{rema}
La condition (i) de la Proposition \ref{nullita} peut s'exprimer par l'égalité suivante:
\begin{align}\label{carattereintegro}
\mathrm{val}_{\Q}(\chi_1(p)) + \mathrm{val}_{\Q}(\chi_2(p)) + |\underline{d}_{S\setminus J}| = 0.
\end{align}
\end{rema}

On termine cette section par quelques remarques sur le cas localement algébrique. Soient $\chi_1,\chi_2 \colon F^{\times} \to E^{\times}$ deux caractères localement constants et  $\underline{d}$ une $|S|$-uplet d'entiers positifs ou nuls. Posons:
\[
I(\chi,\underline{d}) = \Big(  \bigotimes_{\sigma \in S} (\Sym^{d_{\sigma}} E^2)^{\sigma}  \Big) \otimes_E \Big(\Ind_{P}^{G} \chi_1 \otimes \chi_2|\cdot|^{-1} \Big),
\]  
où $\Ind_{P}^{G} \chi_1 \otimes \chi_2|\cdot|^{-1}$ désigne l'induite lisse usuelle. D'après la Proposition \ref{nullita} et d'après \cite[Lemme 7.9]{pas} on connait deux conditions nécessaires pour que $I(\chi,\underline{d})^{\bigwedge}$ soit non nul, c'est-à-dire:
\begin{itemize}
\item[(i)] $\mathrm{val}_{\Q}(\chi_1(p)) + \mathrm{val}_{\Q}(\chi_2(p)) + 1 + |\underline{d}| = 0$; 
\item[(ii)] $\mathrm{val}_{\Q}(\chi_2(p))+ 1 + |\underline{d}| \geq 0$ et  $\mathrm{val}_{\Q}(\chi_1(p)) + 1 + |\underline{d}| \geq 0$.
\end{itemize}
On voit facilement que $(i)$ et $(ii)$ sont équivalents à
\begin{itemize}
\item[(i')] $\mathrm{val}_{\Q}(\chi_1(p)) + \mathrm{val}_{\Q}(\chi_2(p)) + 1 + |\underline{d}| = 0$; 
\item[(ii')] $\mathrm{val}_{\Q}(\chi_2(p)) \leq 0$ et  $\mathrm{val}_{\Q}(\chi_1(p))\leq 0$.
\end{itemize}

Rappelons la conjecture suivante qui est un cas particulier d'une conjecture plus général formulée par  Breuil et Schneider dans \cite{bs}.
\begin{conj}\label{congettura}
Avec les notations précédentes, les conditions $(i')$ et $(ii')$ sont aussi suffisantes pour que $I(\chi,\underline{d})^{\bigwedge}$ soit non nul. 
\end{conj} 

\begin{rema}
{\rm
On connait une réponse positive à la conjecture \ref{congettura} dans les cas suivants:
\begin{itemize} 
\item[$\bullet$] Si $F = \Q$ (\cite[Corollaire 5.3.1]{bb}); 
\item[$\bullet$] Si $\chi_2\chi_1^{-1}$ est un caractère modérément ramifié et $\underline{d} = \underline{0}$  (\cite[Proposition 0.10]{vig}, \cite[Théorème 1.2]{ks}); 
\item[$\bullet$]  Si $\chi_2\chi_1^{-1}$ est un caractère non ramifié et le vecteur d'entiers $\underline{d}$ est sujet à  quelques restrictions (\cite{marco2}). 
\end{itemize}  
}
\end{rema}

\subsection{Passage aux duaux}


Conservons les notations du §\ref{lattiprov} et supposons que les conditions (i) et (ii) de la Proposition \ref{nullita} soient satisfaites ce qui implique en particulier $r \geq 0$. Posons:
\begin{align*}
J' &= J \coprod \{\sigma \in S\backslash J,\, d_{\sigma}+1 > r\},\quad \chi_1' = \chi_1,\quad
\chi_2' = \chi_2 \prod_{\sigma\in J'\setminus J}\sigma^{d_{\sigma}}.
\end{align*}
On a une immersion fermée $G$-équivariante:
\begin{align}\label{canonica}
I(\chi,J,\underline{d}_{S\setminus J}) \into I(\chi',J',\underline{d}_{S\setminus J'}). 
\end{align}
Le résultat suivant donne des indications concernant la structure de $I(\chi,J,\underline{d}_{S\setminus J} )^{\bigwedge}$, ou plus précisément ses vecteurs localement $\Q$-analytiques. 

\begin{prop}\label{dsigma}
Supposons que les conditions de la Proposition \ref{nullita} soient satisfaites. Alors les conditions suivantes sont équivalentes:
\begin{itemize}
\item[(i)] Toute application continue, $E$-linéaire et $G$-équivariante $I(\chi,J,\underline{d}_{S\setminus J}) \to B$, où $B$ est un $G$-Banach unitaire, s'étend de manière unique en une application continue, $E$-linéaire et $G$-équivariante $I(\chi',J',\underline{d}_{S\setminus J'}) \to B$. 
\item[(ii)] L'application canonique $I(\chi,J,\underline{d}_{S\setminus J}) \to I(\chi,J,\underline{d}_{S\setminus J} )^{\bigwedge}$ s'étend de manière unique en une application continue, $E$-linéaire et $G$-équivariante $I(\chi',J',\underline{d}_{S\setminus J'}) \to I(\chi,J,\underline{d}_{S\setminus J} )^{\bigwedge}$.  
\item[(iii)] L'application \eqref{canonica} induit un isomorphisme de $G$-Banach unitaires: 
\[
I(\chi,J,\underline{d}_{S\setminus J} )^{\bigwedge} \stackrel{\sim}{\longrightarrow} I(\chi',J',\underline{d}_{S\setminus J'} )^{\bigwedge}
\]
\end{itemize} 
\begin{proof}
L'équivalence des conditions (i), (ii) et (iii) est clair. Breuil montre (i) sous l'hypothèse supplémentaire que l'application de $I(\chi,J,\underline{d}_{S\setminus J})$ dans $B$ est injective (\cite[Théorème 7.1]{bre}). Une preuve similaire, qui utilise de façon cruciale le \cite[Lemme  6.1]{bre}, permet de démontrer le cas général.   
\end{proof} 
\end{prop}  

D'après la Proposition \ref{dsigma} (iii) donner une description explicite de $I(\chi,J,\underline{d}_{S\setminus J} )^{\bigwedge}$ est équivalente à donner une description explicite de $I(\chi',J',\underline{d}_{S\setminus J'} )^{\bigwedge}$. On peut alors supposer que: 
\begin{align}\label{condizioneaggiuntiva}
\forall \sigma \in S\backslash J, \quad r \geq d_{\sigma}+1
\end{align}
ou ce qui revient au même $J=J'$. 

Rappelons (§\ref{lattiprov}) que le complété unitaire universel de $I(\chi,J,\underline{d}_{S\setminus J})$ est le complété  par rapport au sous-$\Oe[P]$-réseau $\Lambda$, qui est engendré par les vecteurs: 
\begin{align}\label{richiamolatti}
\mathbf{1}_{\OF}(z)z^{\underline{n}_{S\setminus J}}z^{\underline{m}_{J}}, \quad \mathbf{1}_{F-\OF}(z)\chi_2\chi_1^{-1}(z) z^{\underline{d}_{S\setminus J}-\underline{n}_{S\setminus J}}z^{-\underline{m}_{J}} 
\end{align}
pour tout $\underline{0}\leqslant\underline{n}_{S\setminus J} \leqslant \underline{d}_{S\setminus J}$ et tout $\underline{m}_{J} \in \mathbb{Z}_{\geq 0}^{|J|}$. 

Rappelons que $I(\chi, J,\underline{d}_{S\setminus J})^{\vee}$ désigne le dual continu de l'espace $I(\chi, J,\underline{d}_{S\setminus J})$. Si $\mu \in I(\chi, J,\underline{d}_{S\setminus J})^{\vee}$ et $f\in I(\chi, J,\underline{d}_{S\setminus J})$, on note $\int_F f(z)\mu(z)$ l'accouplement et on pose:
\[
\int_U f(z) \mu(z) = \int_F \mathbf{1}_U(z) f(z) \mu(z).
\]
où, si $U$ est un ouvert de $F$, $\mathbf{1}_U$ désigne la fonction caractéristique de $U$.


D'après la Remarque \ref{ossuti} l'application canonique $I(\chi,J,\underline{d}_{S\setminus J} ) \to I(\chi,J,\underline{d}_{S\setminus J})^{\bigwedge}$ est d'image dense. Cela implique que l'on a une injection continue
\begin{align}\label{injcont}
(I(\chi,J,\underline{d}_{S\setminus J} )^{\bigwedge})^{\vee} \into I(\chi,J,\underline{d}_{S\setminus J} )^{\vee}. 
\end{align}
Le résultat suivant donne une caractérisation utile de l'image de l'application \eqref{injcont}.

\begin{prop} \label{funzi} Soit $\mu \in I(\chi,J,\underline{d}_{S\setminus J} )^{\vee}$.  Alors $\mu$ est un élément de $(I(\chi,J,\underline{d}_{S\setminus J})^{\bigwedge})^{\vee}$ si et seulement s'il existe une constante $C_{\mu} \in \mathbb{R}_{\geq 0}$ telle que pour tout $n \in \mathbb{Z}$, tout $a \in F$, tout $\underline{0}\leqslant \underline{n}_{S\setminus J}\leqslant \underline{d}_{S\setminus J}$ et tout $\underline{m}_J \in \mathbb{Z}_{\geq 0}^{|J|}$ on a:
\begin{align}
&\Big|\int_{D(a,n)} (z-a)^{\underline{n}_{S\setminus J}}  (z-a)^{\underline{m}_J} \mu(z)\Big| \leq C_{\mu} q^{n (r-|\underline{n}_{S\setminus J}|-|\underline{m}_J|)}  \label{puno} \\
&\Big|\int_{F\setminus D(a,n+1)} \chi_2\chi_1^{-1}(z-a) (z-a)^{\underline{d}_{S\setminus J}- \underline{n}_{S\setminus J}} (z-a)^{-\underline{m}_J} \mu(z)\Big| \leq C_{\mu} q^{n(|\underline{n}_{S\setminus J}|+|\underline{m}_{J}|-r)}.  \label{pdue}
\end{align}

\begin{proof}

La distribution $\mu$ s'étend en une forme linéaire continue sur  $I(\chi,J,\underline{d}_{S\setminus J})^{\bigwedge}$ si et seulement s'il existe une constante \ $C_{\mu} \in \mathbb{R}_{\geq 0}$ telle que
\begin{align}\label{appart}
\forall f \in \Lambda, \quad    \Big|\int_{F} f(z) \mu(z)\Big| \leq C_{\mu}.
\end{align}
En utilisant \eqref{richiamolatti} et l'identité 
\[
\left[\begin{smallmatrix} {0} & {1} \cr {1} & {0} \end{smallmatrix}\right](\mathbf{1}_{\OF}(z)z^{\underline{n}_{S\setminus J}}z^{\underline{m}_{J}}) = \mathbf{1}_{F-\OF}(z)\chi_2\chi_1^{-1}(z) z^{\underline{d}_{S\setminus J}-\underline{n}_{S\setminus J}}z^{-\underline{m}_{J}}
\]
on déduit immédiatement que \eqref{appart} est équivalente aux deux conditions suivantes: 
\begin{align}
& \big|\mu \big(b  (\mathbf{1}_{\OF}(z)z^{\underline{n}_{S\setminus J}}z^{\underline{m}_{J}})\big)\big| \leq C_{\mu} \label{conduno}\\ 
&\big|\mu \big( b  \left[\begin{smallmatrix} {0} & {1} \cr {1} & {0} \end{smallmatrix}\right] (\mathbf{1}_{\OF}(z)z^{\underline{n}_{S\setminus J}}z^{\underline{m}_{J}})\big)\big| \leq C_{\mu} \label{condue}
\end{align}
pour tout $b \in \Bigl\{\left[\begin{smallmatrix} {\varpi_F^n} & {a} \cr {0} & {1} \end{smallmatrix}\right] \ \mbox{pour} \ n\in \mathbb{Z},a \in F \Bigr\}$, tout  $\underline{0}\leqslant \underline{n}_{S\setminus J}\leqslant \underline{d}_{S\setminus J}$  et tout $\underline{m}_J \in \mathbb{Z}_{\geq 0}^{|J|}$. \\
Or, en appliquant la formule \eqref{azioneanalitica} et d'après \eqref{carattereintegro} on obtient:
\begin{align*}
\Big|\mu\Big(\left[\begin{smallmatrix} {\varpi_F^n} & {a} \cr {0} & {1} \end{smallmatrix}\right] (\mathbf{1}_{\OF}(z) z^{\underline{n}_{S\setminus J}}z^{\underline{m}_{ J}})  \Big)\Big| &= \Big|\mu \Big( \mathbf{1}_{D(a,n)}(z) \chi_2(\varpi_F^n) \varpi_F^{n \underline{d}_{S\setminus J}} \Big(\frac{z-a}{\varpi_F^n}\Big)^{\underline{n}_{S\setminus J}} \Big(\frac{z-a}{\varpi_F^n}\Big)^{\underline{m}_J} \Big)\Big| \\
&=  q^{n (|\underline{n}_{S\setminus J}|+|\underline{m}_J|-r)} \Big|\mu \Big(\mathbf{1}_{D(a,n)}(z) (z-a)^{\underline{n}_{S\setminus J}}(z-a)^{\underline{m}_J}\Big)\Big|
\end{align*}
d'où la condition \eqref{puno}.

Un calcul analogue montre que la condition \eqref{condue} est équivalente à la condition \eqref{pdue}. 
\end{proof}
\end{prop}  


\begin{defi}
On appelle distribution tempérée d'ordre $r$ sur $F$ une forme linéaire continue sur l'espace de Banach $B(\chi,J,\underline{d}_{S\setminus J})$.
\end{defi}

D'après §\ref{ana} on sait que $\mathcal{F}(\OF,J,\underline{d}_{S\setminus J})$ s'injecte de façon continue dans $C^{r}(\OF,J,\underline{d}_{S\setminus J})$ et que l'image de $\mathcal{F}(\OF,J,\underline{d}_{S\setminus J})$ dans $C^{r}(\OF,J,\underline{d}_{S\setminus J})$ est dense. En utilisant le fait que $I(\chi,J,\underline{d}_{S\setminus J})$ (resp. $B(\chi,J,\underline{d}_{S\setminus J})$) s'indentifie topologiquement à deux copies de $\mathcal{F}(\OF,J,\underline{d}_{S\setminus J})$ (resp. $C^{r}(\OF,J,\underline{d}_{S\setminus J})$) on en déduit une injection  $\GL(F)$-équivariante continue:
\[
I(\chi,J,\underline{d}_{S\setminus J}) \into B(\chi,J,\underline{d}_{S\setminus J}),
\]
d'où  a une injection continue:
\begin{align}\label{injcont2}
B(\chi,J,\underline{d}_{S\setminus J})^{\vee} \into I(\chi,J,\underline{d}_{S\setminus J})^{\vee}.
\end{align}

Le résultat suivant donne une caractérisation utile de l'image de l'application \eqref{injcont2}.

\begin{prop} \label{crit} Soit $\mu \in I(\chi,J,\underline{d}_{S\setminus J})^{\vee}$.  Alors $\mu$ est tempérée d'ordre $r$ sur $F$ si et seulement s'il existe une constante $C_{\mu} \in \mathbb{R}_{\geq 0}$ telle que
\begin{align}
&\Big|\int_{D(a,n)} (z-a)^{\underline{n}_{S\setminus J}} (z-a)^{\underline{m}_J} \mu(z)\Big| \leq C_{\mu} q^{n (r-|\underline{n}_{S\setminus J}|-|\underline{m}_J|)}  \label{unocrit} 
\end{align}
pour tout $a \in \varpi_F \OF$, tout $\underline{0}\leqslant \underline{n}_{S\setminus J}\leqslant \underline{d}_{S\setminus J}$ tout $\underline{m}_J \in \mathbb{Z}_{\geq 0}^{|J|}$ et tout $n \geq 1$;
\begin{align} 
&\Big|\int_{F\setminus D(0,n+1)} \chi_2\chi_1^{-1}(z) z^{\underline{d}_{S\setminus J}-\underline{n}_{S\setminus J}}  z^{-\underline{m}_J} \mu(z)\Big| \leq C_{\mu} q^{n (|\underline{n}_{S\setminus J}|+|\underline{m}_J|-r)}  \label{duecrit} 
\end{align}
pour tout $\underline{0} \leqslant \underline{n}_{S\setminus J}\leqslant \underline{d}_{S\setminus J}$, tout $\underline{m}_J \in \mathbb{Z}_{\geq 0}^{|J|}$ et tout $n \leq 0$; 
\begin{align} 
&\Big|\int_{D(\frac{1}{a},n-\frac{2 val_F(a)}{f})} \chi_2\chi_1^{-1}(z) z^{ \underline{d}_{S\setminus J}} \Big(\frac{1}{z}-a\Big)^{\underline{n}_{S\setminus J}} \Big(\frac{1}{z}-a\Big)^{\underline{m}_J} \mu(z) \Big| \leq C_{\mu} q^{n (r-|\underline{n}_{S\setminus J}|-|\underline{m}_J|)} \label{trecrit} 
\end{align}
pour tout $a \in \OF-\{0 \}$, tout $\underline{0} \leqslant \underline{n}_{S\setminus J}\leqslant \underline{d}_{S\setminus J}$, tout $\underline{m}_J \in \mathbb{Z}_{\geq 0}^{|J|}$ et tout entier $n> \frac{val_F(a)}{f}$. 
\begin{proof}

L'application \eqref{traduiso} (resp. \eqref{isocr}) induit un isomorphisme topologique de $I(\chi, J,\underline{d}_{S\setminus J})^{\vee}$ dans $(\mathcal{F}(\OF, J,\underline{d}_{S\setminus J})^{\vee})^2$ (resp. de $B(\chi, J,\underline{d}_{S\setminus J})^{\vee}$ dans $(C^r(\OF, J,\underline{d}_{S\setminus J})^{\vee})^2$). Si l'on note $(\mu_1,\mu_2)$ l'élément de $(\mathcal{F}(\OF, J,\underline{d}_{S\setminus J})^{\vee})^2$ qui correspond à $\mu$ via cet isomorphisme alors il est clair que $\mu$ est tempérée d'ordre $r$ sur $F$ si et seulement si les distributions $\mu_1$ et $\mu_2$ sont tempérées d'ordre $r$ sur $\OF$. D'après le Théorème \ref{velu}, la distribution $\mu_1$ (resp.  $\mu_2$) est tempérée d'ordre $r$ sur $\OF$ si et seulement s'il existe une constante $C_{\mu_1} \in \mathbb{R}_{\geq 0}$ (resp. $C_{\mu_2} \in \mathbb{R}_{\geq 0}$) telle que pour tout $a \in \OF$, tout $\underline{0}\leqslant \underline{n}_{S\setminus J} \leqslant \underline{d}_{S\setminus J}$, tout $\underline{m}_J \in \mathbb{Z}_{\geq 0}^{|J|}$ et tout $n \geq 0$ on a:
\begin{align}
\Big| \mu_1 \Big(\mathbf{1}_{D(a,n)}(z) (z-a)^{\underline{n}_{S\setminus J}} (z-a)^{\underline{m}_J}  \Big)\Big| &\leq C_{\mu_1} q^{n (r- |\underline{n}_{S\setminus J}|-|\underline{m}_J|)} \label{unounouno} \\
\Big| \mu_2 \Big(\mathbf{1}_{D(a,n)}(z) (z-a)^{\underline{n}_{S\setminus J}} (z-a)^{\underline{m}_J}  \Big)\Big| &\leq C_{\mu_2} q^{n (r- |\underline{n}_{S\setminus J}|-|\underline{m}_J|)}. \label{dueduedue}  
\end{align}

La fonction $f$ correspondant via \eqref{traduiso2} au couple 
\[
(f_1,f_2) = (\mathbf{1}_{D(a,n)}(z) (z-a)^{\underline{n}_{S\setminus J}} (z-a)^{\underline{m}_J},0) 
\]     
est la fonction $\mathbf{1}_{D(\varpi_F a,n+1)}(z) \big(\frac{z}{\varpi_F}-a\big)^{\underline{n}_{S\setminus J}} \big(\frac{z}{\varpi_F}-a \big)^{\underline{m}_J}$ et donc la condition \eqref{unounouno} se traduit par
\[
\Big| \mu\Big( \mathbf{1}_{D(\varpi_F a,n+1)}(z) (z-\varpi_F a)^{\underline{n}_{S\setminus J}} (z- \varpi_F a)^{\underline{m}_J}   \Big)  \Big| \leq C_{\mu_1} q^{(n+1) (r- |\underline{n}_{S\setminus J}|-|\underline{m}_J|)}
\]
pour tout $a \in \OF$, tout $\underline{0}\leqslant \underline{n}_{S\setminus J} \leqslant \underline{d}_{S\setminus J}$, tout $\underline{m}_J \in \mathbb{Z}_{\geq 0}^{|J|}$ et tout $n \geq 0$,  d'où \eqref{unocrit}. 

La fonction $f$ correspondant via \eqref{traduiso2} au couple 
\[
(f_1,f_2) = (0,\mathbf{1}_{D(a,n)}(z) (z-a)^{\underline{n}_{S\setminus J}} (z-a)^{\underline{m}_J}) 
\]     
est la fonction $\mathbf{1}_{\{z: \ |\frac{1}{z}-a| \leq |\varpi_F^n| \}}(z) \chi_2\chi_1^{-1}(z) z^{\underline{d}_{S\setminus J}} \big(\frac{1}{z}-a\big)^{\underline{n}_{S\setminus J}} \big(\frac{1}{z}-a \big)^{\underline{m}_J}$. On va distinguer deux cas. 
\begin{itemize}
\item[$\bullet$] Si $a \in D(0,n)$ on a $\{z: \ |\frac{1}{z}-a| \leq |\varpi_F^n|  \} = F\backslash D(0,-n+1)$ et donc la condition \eqref{dueduedue} se traduit par
\begin{align}\label{tradi}
\Big| \mu\Big( \mathbf{1}_{F\setminus D(0,n+1)}(z) \chi_2\chi_1^{-1}(z) z^{\underline{d}_{S\setminus J}} \Big(\frac{1}{z}-a\Big)^{\underline{n}_{S\setminus J}} \Big(\frac{1}{z}-a \Big)^{\underline{m}_J}\Big)  \Big| \leq C_{\mu_2} q^{n ( |\underline{n}_{S\setminus J}|+|\underline{m}_J|-r)}
\end{align}
pour tout $\underline{0} \leqslant \underline{n}_{S\setminus J}\leqslant \underline{d}_{S\setminus J}$, tout $\underline{m}_J \in \mathbb{Z}_{\geq 0}^{|J|}$ et tout $n \leq 0$. En développant $\big(\frac{1}{z}-a\big)^{\underline{n}_{S\setminus J}}$ et $ \big(\frac{1}{z}-a \big)^{\underline{m}_J}$ on voit facilement que la condition \eqref{tradi} est équivalente à la condition \eqref{duecrit}.  
\item[$\bullet$] Si $a \in \OF \backslash D(0,n)$ on a $\{z: \ |\frac{1}{z}-a| \leq |\varpi_F^n|  \} = D(\frac{1}{a},n-\frac{2 val_F(a)}{f})$ et la condition \eqref{dueduedue} se traduit par la condition \eqref{trecrit}. 
\end{itemize}
\end{proof}
\end{prop}

\begin{coro}\label{corollariofacile}
Soit $\mu \in I(\chi,J,\underline{d}_{S\setminus J})^{\vee}$ . Alors $\mu$ est dans $\Pi(\chi,J,\underline{d}_{S\setminus J})^{\vee}$ si et seulement s'il existe une constante $C_{\mu} \in \mathbb{R}_{\geq 0}$ vérifiant \eqref{unocrit}, \eqref{duecrit}, \eqref{trecrit} et les deux conditions supplémentaires suivantes:  
\begin{align}
&\int_{F} z^{\underline{n}_{S\setminus J}} z^{\underline{m}_J}\mu(z) = 0; \label{funo} \\
&\int_{F} \chi_2\chi_1^{-1}(z-a)(z-a)^{\underline{d}_{S\setminus J}- \underline{n}_{S\setminus J}} (z-a)^{-\underline{m}_J} \mu(z) = 0 \label{fdue}
\end{align}
pour tout $a\in F$, tout $\underline{0}\leqslant \underline{n}_{S\setminus J}\leqslant \underline{d}_{S\setminus J}$ et tout $\underline{m}_J \in \mathbb{Z}_{\geq 0}^{|J|}$ tels que $r-(|\underline{n}_{S\setminus J}|+|\underline{m}_J|)>0$.
\begin{proof}
C'est une conséquence immédiate de la proposition \ref{crit} et du Lemme \ref{funzionicr}.
\end{proof}
\end{coro} 
 
\section{Preuve du Théorème principal}

Conservons les notations du §\ref{lattiprov} et supposons que les conditions (i) et (ii) de la Proposition \ref{nullita} et la condition \eqref{condizioneaggiuntiva} soient satisfaites. Nous nous proposons de montrer que les conditions \eqref{puno} et \eqref{pdue} sélectionnent exactement les distributions tempérées d'ordre $r$ sur $F$ annulant les fonctions $z^{\underline{n}_{S\setminus J}} z^{\underline{m}_{J}}$ et $\chi_2\chi_1^{-1}(z-a) (z-a)^{\underline{d}_{S\setminus J}- \underline{n}_{S\setminus J}} (z-a)^{-\underline{m}_{J}}$ pour tout $a \in F$, tout $\underline{m}_J \in \mathbb{Z}_{\geq 0}^{|J|}$  et tout $\underline{0}\leqslant \underline{n}_{S\setminus J} \leqslant \underline{d}_{S\setminus J}$ tels que $r - (|\underline{n}_{S\setminus J}|+ |\underline{m}_J|) > 0$. Plus précisément:

\begin{theo}\label{teoremaprincipale}
Soit $\mu \in I(\chi,J,\underline{d}_{S\setminus J})^{\vee}$. Les deux conditions suivantes sont équivalentes:
\begin{itemize}
\item[(A)] La distribution $\mu$ vérifie les conditions \eqref{puno} et \eqref{pdue};
\item[(B)] La distribution $\mu$ vérifie les conditions \eqref{unocrit}, \eqref{duecrit}, \eqref{trecrit}, \eqref{funo} et \eqref{fdue}. 
\end{itemize}
\end{theo}

\subsection{(A) $\Rightarrow $ (B)}

Supposons que $\mu$ vérifie \eqref{puno} et \eqref{pdue}. Alors \textit{a fortiori} $\mu$ vérifie \eqref{unocrit} et \eqref{duecrit}. Montrons que \eqref{puno} implique \eqref{trecrit} quitte à changer $C_{\mu}$. Pour cela on aura besoin de l'équivalence suivante.  


\begin{lemm}\label{unadelle1}
Quitte à modifier la constante $C_{\mu}$ les deux conditions suivantes sont équivalentes:
\begin{itemize}
\item[(i)] La condition \eqref{trecrit};
\item[(ii)] Il existe un entier $n_0>0$ tel que \eqref{trecrit} est satisfaite pour tout $a \in \OF-\{0\}$, tout $\underline{0} \leqslant \underline{n}_{S\setminus J}\leqslant \underline{d}_{S\setminus J}$, tout $\underline{m}_J \in \mathbb{Z}_{\geq 0}^{|J|}$ et tout $n > n_0+\frac{val_F(a)}{f}$. 
\end{itemize}
\begin{proof}
$(i) \Rightarrow (ii)$ est immédiat. 

Montrons $(ii) \Rightarrow (i)$. Soit $a\in \OF-\{0\}$ et $ \frac{val_F(a)}{f} < n \leq n_0+ \frac{val_F(a)}{f}$. Si l'on note $n' = n+n_0$ on peut écrire $D\big(\frac{1}{a},n-\frac{2 val_F(a)}{f}\big)$ comme union de disques de la forme $D' = D\big(\frac{1}{a'},n'-\frac{2 val_F(a)}{f}\big)$ avec $|a| = |a'|$ (et donc $|a-a'| \leq q^{-n}$). En écrivant $\big(\frac{1}{z}-a\big)^{\underline{i}} = \big( \big(\frac{1}{z}-a'\big)+\big(a'-a\big)\big)^{\underline{i}}$ avec $ \underline{i} \in \{\underline{n}_{S\setminus J},\underline{m}_J \} $  et en développant on obtient:
\begin{align*}
&\quad \quad \Big|\mu\Big(\mathbf{1}_{D'}(z) \chi_2\chi_1^{-1}(z) z^{ \underline{d}_{S\setminus J}} \Big(\frac{1}{z}-a\Big)^{\underline{n}_{S\setminus J}} \Big(\frac{1}{z}-a\Big)^{\underline{m}_J} \Big) \Big| \\
\stackrel{\phantom{(ii)}}\leq & \  \sup_{\substack{
\underline{0}\leqslant \underline{k}_{S\setminus J} \leqslant \underline{n}_{S\setminus J} \\
\underline{0}\leqslant \underline{l}_{J} \leqslant \underline{m}_{J}}} \Big\{      |a-a'|^{|\underline{n}_{S\setminus J}|-|\underline{k}_{S\setminus J}|+|\underline{m}_{J}|-|\underline{l}_{J}|} \\
& \qquad \qquad \quad \qquad \qquad \qquad \qquad \qquad \cdot \Big|\mu\Big(\mathbf{1}_{D'}(z) \chi_2\chi_1^{-1}(z) z^{ \underline{d}_{S\setminus J}} \Big(\frac{1}{z}-a'\Big)^{\underline{k}_{S\setminus J}} \Big(\frac{1}{z}-a'\Big)^{\underline{l}_J} \Big) \Big| \Big\} \\
 \stackrel{(ii)}\leq & \ \sup_{\substack{
\underline{0}\leqslant \underline{k}_{S\setminus J} \leqslant \underline{n}_{S\setminus J} \\
\underline{0}\leqslant \underline{l}_{J} \leqslant \underline{m}_{J}}} q^{n(-|\underline{n}_{S\setminus J}|+|\underline{k}_{S\setminus J}|-|\underline{m}_{J}|+|\underline{l}_{J}|)} C_{\mu} q^{n' (r-|\underline{k}_{S\setminus J}|-|\underline{l}_J|)} \\
\stackrel{\phantom{(ii)}}= & \ C_{\mu} q^{n(r-|\underline{n}_{S\setminus J}|-|\underline{m}_{ J}|)}  q^{(n'-n) r} \\
\stackrel{\phantom{(ii)}}\leq & \ C_{\mu}' q^{n(r-|\underline{n}_{S\setminus J}|-|\underline{m}_{ J}|)}
\end{align*}
où l'on a posé $C_{\mu}'\ugu C_{\mu} q^{n_0 r}$. Comme le dernier terme de dépend pas du choix de $a$ on peut conclure.

\end{proof}
\end{lemm}

\begin{prop}\label{messi1}
Quitte à modifier la constante $C_{\mu}$ la condition \eqref{puno} implique la condition \eqref{trecrit}.
\begin{proof}
Notons $n_0$ le plus petit entier positif tel que $(\chi_2 \chi_1^{-1})|_{D(1,n_0)}$ est une fonction $J$-analytique. D'après le Lemme \ref{unadelle1} il suffit de montrer que la condition \eqref{trecrit} est satisfaite  pour tout $a \in \OF-\{0\}$, tout $\underline{0} \leqslant \underline{n}_{S\setminus J}\leqslant \underline{d}_{S\setminus J}$, tout $\underline{m}_J \in \mathbb{Z}_{\geq 0}^{|J|}$ et tout $n > n_0+\frac{val_F(a)}{f}$. Posons $D = D\big(\frac{1}{a},n-\frac{2 val_F(a)}{f}\big)$. 

D'après l'égalité: 
\[
\mathbf{1}_D(z) \Big(\frac{1}{z}-a  \Big)^{\underline{n}_{S\setminus J}} = \mathbf{1}_D(z) (-1)^{\underline{n}_{S\setminus J}}  z^{-\underline{n}_{S\setminus J}} a^{\underline{n}_{S\setminus J}} \Big(z- \frac{1}{a}  \Big)^{\underline{n}_{S\setminus J}} 
\]
et, en écrivant $z^{\underline{d}_{S\setminus J}-\underline{n}_{S\setminus J}}  =  (z-\frac{1}{a}+\frac{1}{a})^{\underline{d}_{S\setminus J}-\underline{n}_{S\setminus J}}$ et en développant on obtient:
\[
\mathbf{1}_D(z) z^{\underline{d}_{S\setminus J}} \Big(\frac{1}{z}-a  \Big)^{\underline{n}_{S\setminus J}} = \mathbf{1}_D(z) \sum_{\underline{0}\leqslant \underline{k}_{S\setminus J} \leqslant \underline{d}_{S\setminus J}- \underline{n}_{S\setminus J}} \mu_{\underline{k}_{S\setminus J}} a^{-\underline{k}_{S\setminus J}+\underline{n}_{S\setminus J}}\Big(z- \frac{1}{a}  \Big)^{\underline{d}_{S\setminus J}-\underline{k}_{S\setminus J}},
\]
où les $\mu_{\underline{k}_{S\setminus J}}$ sont des entiers. De manière analogue, en écrivant $z^{-\underline{m}_J} = (z-\frac{1}{a}+\frac{1}{a})^{-\underline{m}_J}$  et en développant on a:
\[
\mathbf{1}_D(z) z^{-\underline{m}_J} = \mathbf{1}_D(z) a^{\underline{m}_J}  \sum_{\underline{r}_J \geqslant \underline{0}} \lambda_{\underline{r}_J} a^{\underline{r}_J} \Big(z-\frac{1}{a}  \Big)^{\underline{r}_J}
\]
où les $\lambda_{\underline{r}_J}$ sont des entiers, d'où l'égalité
\begin{align*}
\mathbf{1}_D(z) \Big(\frac{1}{z}-a  \Big)^{\underline{m}_J} &= \mathbf{1}_D(z) (-1)^{\underline{m}_J} z^{-\underline{m}_J}a^{\underline{m}_J} \Big(z-\frac{1}{a}  \Big)^{\underline{m}_J} \\
&= \mathbf{1}_D(z)   \sum_{\underline{r}_J \geqslant \underline{0}} \lambda_{\underline{r}_J} a^{2 \underline{m}_J+\underline{r}_J} \Big(z-\frac{1}{a}  \Big)^{\underline{m}_J+\underline{r}_J}.
\end{align*}
Remarquons que
\[
z \in D \Rightarrow az \in D\Big(1,n-\frac{val_F(a)}{f}\Big) \subseteq D(1,n_0)
\]
ce qui implique 
\begin{align*}
\mathbf{1}_D(z) \chi_2 \chi_1^{-1}(z) &= \chi_2 \chi_1^{-1}(a^{-1})\mathbf{1}_D(z) \chi_2 \chi_1^{-1}(az) \\
&=  \chi_2 \chi_1^{-1}(a^{-1})\mathbf{1}_D(z) \sum_{\underline{l}_J \geqslant \underline{0}} b_{\underline{l}_J}(az-1)^{\underline{l}_J} \\
&= \chi_2 \chi_1^{-1}(a^{-1})\mathbf{1}_D(z) \sum_{\underline{l}_J \geqslant \underline{0}} b_{\underline{l}_J} a^{\underline{l}_J} \Big(z-\frac{1}{a}\Big)^{\underline{l}_J}.
\end{align*}
avec $b_{\underline{l}_J} \in E$ et $|b_{\underline{l}_J}| q^{-n_0} \to 0$ quand $|\underline{l}_J| \to +\infty$. Notons $C = \sup_{\underline{l}_J} |b_{\underline{l}_J}|$ et remarquons que d'après \eqref{carattereintegro} on a $|\chi_2 \chi_1^{-1}(a^{-1})| = |a|^{|\underline{d}_{S\setminus J}|-2r}$.

Par les égalités précédentes on obtient:
\begin{align*}
& \ \Big|\mu\Big(\mathbf{1}_D(z) \chi_2\chi_1^{-1}(z) z^{ \underline{d}_{S\setminus J}} \Big(\frac{1}{z}-a\Big)^{\underline{n}_{S\setminus J}} \Big(\frac{1}{z}-a\Big)^{\underline{m}_J} \Big)\Big| \\
\leq & \ C |a|^{|\underline{d}_{S\setminus J}|-2r} \sup_{\substack{
\underline{0} \leqslant \underline{k}_{S\setminus J} \leqslant  \underline{d}_{S\setminus J}-\underline{n}_{S\setminus J}  \\
\underline{l}_{J} \geqslant \underline{0},\ \underline{r}_{J} \geqslant \underline{0}
}}  \Big\{ |a|^{2|\underline{m}_J|+|\underline{r}_J|+|\underline{l}_J|-|\underline{k}_{S\setminus J}|+|\underline{n}_{S\setminus J}|}   \\ 
& \qquad \qquad \qquad \qquad \qquad \qquad \qquad \qquad \quad    \cdot  \Big|\mu\Big( \mathbf{1}_D(z) \Big(z-\frac{1}{a}  \Big)^{\underline{d}_{S\setminus J}-\underline{k}_{S\setminus J}} \Big(z-\frac{1}{a}  \Big)^{\underline{m}_J+\underline{l}_J+\underline{r}_J}   \Big)\Big| \Big\} 
\end{align*}
et comme \eqref{puno} implique l'inégalité
\[
\Big|\mu\Big( \mathbf{1}_D(z) \Big(z-\frac{1}{a}  \Big)^{\underline{d}_{S\setminus J}-\underline{k}_{S\setminus J}} \Big(z-\frac{1}{a}  \Big)^{\underline{m}_J+\underline{l}_J+\underline{r}_J}   \Big)\Big| \leq C_{\mu} \Big|\frac{\varpi_F^n}{a^2} \Big|^{|\underline{d}_{S\setminus J}|-|\underline{k}_{S\setminus J}|+|\underline{m}_J|+|\underline{l}_J|+|\underline{r}_J|-r}
\]
on en déduit
\[
\Big|\mu\Big(\mathbf{1}_D(z) \chi_2\chi_1^{-1}(z) z^{ \underline{d}_{S\setminus J}} \Big(\frac{1}{z}-a\Big)^{\underline{n}_{S\setminus J}} \Big(\frac{1}{z}-a\Big)^{\underline{m}_J} \Big)\Big| \leq C C_{\mu} q^{n(r-|\underline{n}_{S\setminus J}|-|\underline{m}_J|)},
\]
d'où le résultat.
\end{proof}
\end{prop}

D'après la Proposition \ref{messi1} on peut étendre $\mu$ en une distribution tempérée d'ordre $r$ sur $F$. Il reste à montrer que $\mu$, vu dans   $B(\chi,J,\underline{d}_{S\setminus J})^{\vee}$, annule l'espace  $L(\chi,J,\underline{d}_{S\setminus J})$. Or, d'après \eqref{puno} on a pour tout $\underline{0}\leqslant \underline{n}_{S\setminus J}\leqslant \underline{d}_{S\setminus J}$ et tout $\underline{m}_J \in \mathbb{Z}_{\geq 0}^{|J|}$ tels que $r-(|\underline{n}_{S\setminus J}|+|\underline{m}_J|)>0$:
\[
\Big|\int_{D(0,n)} z^{\underline{n}_{S\setminus J}} z^{\underline{m}_J} \mu(z)\Big|  \to 0 \quad \mbox{quand} \quad n\to -\infty 
\]  
et d'après \eqref{pdue} on a pour tout $a\in F$, tout $\underline{0}\leqslant \underline{n}_{S\setminus J}\leqslant \underline{d}_{S\setminus J}$ et tout $\underline{m}_J \in \mathbb{Z}_{\geq 0}^{|J|}$ tels que $r-(|\underline{n}_{S\setminus J}|+|\underline{m}_J|)>0$:
\[
\Big|\int_{F\setminus D(a,n+1)} \chi_2\chi_1^{-1}(z-a) (z-a)^{\underline{d}_{S\setminus J}- \underline{n}_{S\setminus J}} (z-a)^{-\underline{m}_J} \mu(z)\Big|  \to 0 \quad \mbox{quand} \quad n\to +\infty, 
\]
d'où le résultat, qui permet de terminer la preuve de $(A) \Rightarrow (B)$.

\subsection{(B) $\Rightarrow $ (A)}

Montrer que les conditions \eqref{unocrit}, \eqref{duecrit}, \eqref{trecrit}, \eqref{funo} et  \eqref{fdue} impliquent les conditions \eqref{puno} et \eqref{pdue} requiert quelques préliminaires. Commençons par donner une description équivalente des conditions \eqref{puno} et \eqref{pdue}.

\begin{lemm}\label{condisupp1}
La condition \eqref{puno} est satisfaite (quitte à changer $C_{\mu}$) si et seulement si les trois conditions suivantes sont vérifiées: 
\begin{itemize}
\item[(i)] \eqref{puno} pour tout $a \in F$ et tout $n \in \mathbb{Z}$ tels que  $D(a,n)\cap \varpi_F \OF = \emptyset$, tout $0\leqslant \underline{n}_{S\setminus J}\leqslant \underline{d}_{S\setminus J}$ et tout $\underline{m}_J \in \mathbb{Z}_{\geq 0}^{|J|}$;
\item[(ii)] \eqref{puno} pour tout $a \in \varpi_F\OF$, tout $n \in \mathbb{Z}_{\geq 1}$, tout $0\leqslant \underline{n}_{S\setminus J}\leqslant \underline{d}_{S\setminus J}$ et tout $\underline{m}_J \in \mathbb{Z}_{\geq 0}^{|J|}$;
\item[(iii)] \eqref{puno} pour $a=0$, pour tout entier $n \leq 0$, tout $0\leqslant \underline{n}_{S\setminus J}\leqslant \underline{d}_{S\setminus J}$ et tout $\underline{m}_J \in \mathbb{Z}_{\geq 0}^{|J|}$  tels que $r- (|\underline{n}_{S\setminus J}|+|\underline{m}_J|)>0$.
\end{itemize}
\begin{proof}

$\eqref{puno} \Rightarrow (i),(ii), (iii)$ est immédiat.

Montrons $(i),(ii), (iii) \Rightarrow \eqref{puno}$. Il suffit de vérifier la condition \eqref{puno} pour $a=0$, pour tout entier $n \leq 0$, tout $0\leqslant \underline{n}_{S\setminus J}\leqslant \underline{d}_{S\setminus J}$ et tout $\underline{m}_J \in \mathbb{Z}_{\geq 0}^{|J|}$  tels que $r- (|\underline{n}_{S\setminus J}|+|\underline{m}_J|)\leq 0$. 

Notons $R\subset \OF$ un système de représentants des classes de $\OF/\varpi_F \OF$ tel que $0 \in R$ et fixons $m \in \mathbb{Z}_{>0}$ tel que $n+m > 0$. Donc on a:
\begin{align*}
\mathbf{1}_{D(0,n)}(z) z^{\underline{n}_{S\setminus J}} z^{\underline{m}_{J}} = \mathbf{1}_{D(0,n+m)}(z) z^{\underline{n}_{S\setminus J}} z^{\underline{m}_{J}} + \sum_{j=0}^{m-1} \sum_{a_i \in R-\{0\}}\mathbf{1}_{D(a_i \varpi_F^{n+j}, n+j+1)}(z) z^{\underline{n}_{S\setminus J}} z^{\underline{m}_{J}}. 
\end{align*}
En utilisant $(ii)$ et $r-(|\underline{n}_{S\setminus J}|+|\underline{m}_J|) \leq 0$ on obtient:
\[
\Big|\mu\Big(\mathbf{1}_{D(0,n+m)}(z) z^{\underline{n}_{S\setminus J}} z^{\underline{m}_{J}}  \Big)\Big| \leq C_{\mu}q^{(n+m)(r-|\underline{n}_{S\setminus J}|-|\underline{m}_J|)} \leq C_{\mu}q^{n(r-|\underline{n}_{S\setminus J}|-|\underline{m}_J|)}.
\]
Il reste à minorer les termes de la somme. Soit $a_i \in R-\{0\}$ et $0\leq j \leq m-1$. En écrivant $z^{\underline{n}_{S\setminus J}} = (z-a_i\varpi_F^{n+j}+a_i\varpi_F^{n+j})^{\underline{n}_{S\setminus J}}$ (resp. $z^{\underline{m}_{J}} = (z-a_i\varpi_F^{n+j}+a_i\varpi_F^{n+j})^{\underline{m}_{J}}$) et en développant on obtient:
\begin{align*}
&\quad \ \ \Big|\mu\Big(\mathbf{1}_{D(a_i \varpi_F^{n+j}, n+j+1)}(z) z^{\underline{n}_{S\setminus J}} z^{\underline{m}_{J}}\Big)\Big| \\
\stackrel{\phantom{(i)}}\leq & \sup_{\substack{\underline{0} \leqslant \underline{l}_{S\setminus J}\leqslant \underline{n}_{S\setminus J}  \\
\underline{0} \leqslant \underline{k}_J \leqslant \underline{m}_J}} \Big\{ \Big| \mu\Big(\mathbf{1}_{D(a_i  \varpi_F^{n+j},n+j+1)}(z) (a_i\varpi_F^{n+j})^{\underline{l}_{S\setminus J}} (a_i\varpi_F^{n+j})^{\underline{k}_J} \\ 
& \qquad \qquad \qquad \qquad \qquad\qquad \qquad \qquad \quad  \cdot (z-a_i\varpi_F^{n+j})^{\underline{n}_{S\setminus J}-\underline{l}_{S\setminus J}} (z-a_i\varpi_F^{n+j})^{\underline{m}_{ J}-\underline{k}_J} \Big) \Big| \Big\} \\
\stackrel{(i)}\leq & \sup_{\substack{\underline{0} \leqslant \underline{l}_{S\setminus J}\leqslant \underline{n}_{S\setminus J}  \\
\underline{0} \leqslant \underline{k}_J \leqslant \underline{m}_J}}      q^{-(n+j)(|\underline{l}_{S\setminus J}|+|\underline{k}_J|)} C_{\mu} q^{(n+j+1)(r - |\underline{n}_{S\setminus J}|+ |\underline{l}_{S\setminus J}|- |\underline{m}_{J}| + |\underline{k}_J|)} \\
\stackrel{\phantom{(i)}} \leq & \  C_{\mu} q^r q^{(n+j)(r - |\underline{n}_{S\setminus J}|- |\underline{m}_{J}| )}. 
\end{align*}
Comme $r- (|\underline{n}_{S\setminus J}|+|\underline{m}_J|) \leq 0$ on a: 
\[
q^{(n+j)(r - |\underline{n}_{S\setminus J}|- |\underline{m}_{J}|)} \leq q^{n(r - |\underline{n}_{S\setminus J}|- |\underline{m}_{J}|)},
\]
d'où le résultat.
\end{proof}
\end{lemm}

Rappelons que pour tout $k \in \mathbb{Z}_{>0}$ on désigne par $S_k \subset \OF^{\times}$ un système de représentants des classes de $(\OF/\varpi_F^{k}\OF)^{\times}$ et que $l$ désigne le plus petit entier positif tel que ${\chi_1}|_{D(a_i,l)}$ (resp. \ ${\chi_2}|_{D(a_i,l)}$) 
est une fonction $J$-analytique sur l'ouvert $D(a_i,l)$ pour tout $a_i \in S_{l}$. Notons $D(a,n,n+1) = D(a,n)\backslash D(a,n+1)$ pour tout $a \in F$ et tout $n \in \mathbb{Z}$.

\begin{lemm}\label{condisupp2}
Supposons que la condition \eqref{puno} soit satisfaite. Alors la condition \eqref{pdue} est satisfaite si et seulement si les deux conditions suivantes sont vérifiées: 
\begin{itemize}
\item[(i)] \eqref{pdue} pour tout $a \in F$, tout $n\geq 0$, tout $\underline{0}\leqslant \underline{n}_{S\setminus J} \leqslant \underline{d}_{S\setminus J}$ et tout $\underline{m}_J \in \mathbb{Z}_{\geq 0}^{|J|}$ tels que $r-(|\underline{n}_{S\setminus J}|+|\underline{m}_J|) > 0$;
\item[(ii)] \eqref{pdue} pour $a = 0$, pour tout $n\leq 0$, tout $\underline{0}\leqslant \underline{n}_{S\setminus J} \leqslant \underline{d}_{S\setminus J}$ et tout $\underline{m}_J \in \mathbb{Z}_{\geq 0}^{|J|}$ tels que $r-(|\underline{n}_{S\setminus J}|+|\underline{m}_J|) \leq 0$.
\end{itemize}
\begin{proof}
$\eqref{pdue} \Rightarrow (i),(ii)$ est immédiat.

Montrons $(i),(ii) \Rightarrow \eqref{pdue}$. Il suffit de vérifier la condition \eqref{pdue} dans les cas suivants:
\begin{itemize}
\item[$\bullet$]  pour tout $a \in F$, tout $n< 0$, tout $\underline{0}\leqslant \underline{n}_{S\setminus J} \leqslant \underline{d}_{S\setminus J}$ et tout $\underline{m}_J \in \mathbb{Z}_{\geq 0}^{|J|}$ tels que $r-(|\underline{n}_{S\setminus J}|+|\underline{m}_J|) > 0$;
\item[$\bullet$] pour tout $a \neq 0$, tout $n\in \mathbb{Z}$, tout $\underline{0}\leqslant \underline{n}_{S\setminus J} \leqslant \underline{d}_{S\setminus J}$ et tout $\underline{m}_J \in \mathbb{Z}_{\geq 0}^{|J|}$ tels que $r-(|\underline{n}_{S\setminus J}|+|\underline{m}_J|) \leq 0$;
\item[$\bullet$] pour $a = 0$, pour tout $n > 0$, tout $\underline{0}\leqslant \underline{n}_{S\setminus J} \leqslant \underline{d}_{S\setminus J}$ et tout $\underline{m}_J \in \mathbb{Z}_{\geq 0}^{|J|}$ tels que $r-(|\underline{n}_{S\setminus J}|+|\underline{m}_J|) \leq 0$.
\end{itemize}

Remarquons d'abord qu'en utilisant l'égalité:
\[
\forall a\in F, n\in \mathbb{Z}, \quad \mathbf{1}_{D(a,n,n+1)}  = \sum_{a_i \in S_l} \mathbf{1}_{D(a+a_i \varpi_F^{n},n+l)}
\]
un raisonnement analogue à celui du lemme \ref{funzionicr} permet de montrer, en utilisant \eqref{puno}, que pour tout $a \in F$, tout $n\in \mathbb{Z}$, tout $\underline{0}\leqslant \underline{n}_{S\setminus J} \leqslant \underline{d}_{S\setminus J}$ et tout $\underline{m}_J \in \mathbb{Z}_{\geq 0}^{|J|}$ on a:
\begin{align}\label{inelunga}
\Big|\mu\Big(\mathbf{1}_{D(a,n,n+1)}(z) \chi_2\chi_1^{-1}(z-a) (z-a)^{\underline{d}_{S\setminus J}-\underline{n}_{S\setminus J}} (z-a)^{-\underline{m}_J}   \Big)\Big| \leq C_{\mu} q^{n(|\underline{n}_{S\setminus J}|+ |\underline{m}_{J}|-r)},
\end{align}
quitte à modifier $C_{\mu}$. 

\underline{\textit{Premier cas}}. Soit $n< 0$ et fixons $m \in \mathbb{Z}_{\geq 1}$ de sorte que que $n+m > 0$. En utilisant l'égalité
\[
\forall a\in F, \quad \mathbf{1}_{F \backslash D(a,n)} = \mathbf{1}_{F \backslash D(a,n+m)} -  \sum_{j=0}^{m-1} \mathbf{1}_{D(a,n+j,n+j+1)}
\]
on déduit le premier cas de (i) et de \eqref{inelunga}.  

\underline{\textit{Deuxième cas}}. Soit $a\neq 0$ et $n\in \mathbb{Z}$. Choisissons $m\in \mathbb{Z}$ de sorte que $n-m < 0$ et $F \backslash D(a,n-m) = F \backslash D(0,n-m)$. En utilisant l'égalité
\[
\mathbf{1}_{F \backslash D(a,n)} = \mathbf{1}_{F \backslash D(a,n-m)} +  \sum_{j=0}^{m+1} \mathbf{1}_{D(a,n-m-j,n-m-j+1)}
\]
on déduit le deuxième cas de (ii) et de \eqref{inelunga}.

\underline{\textit{Troisième cas}}. Le même raisonnement que pour le deuxième cas s'applique.

\end{proof}
\end{lemm}

Remarquons que \eqref{unocrit} est exactement \eqref{puno} pour tout $a \in \varpi_F\OF$, tout $n \in \mathbb{Z}_{\geq 1}$, tout $0\leqslant \underline{n}_{S\setminus J}\leqslant \underline{d}_{S\setminus J}$ et tout $\underline{m}_J \in \mathbb{Z}_{\geq 0}^{|J|}$ et que \eqref{duecrit} est exactement \eqref{pdue} pour $a = 0$, pour tout $n\leq 0$, tout $\underline{0}\leqslant \underline{n}_{S\setminus J} \leqslant \underline{d}_{S\setminus J}$ et tout $\underline{m}_J \in \mathbb{Z}_{\geq 0}^{|J|}$. D'après les Lemmes \ref{condisupp1} et \ref{condisupp2} il reste alors à montrer:
\begin{itemize}
\item[(i)] \eqref{puno} pour tout $a \in F$ et tout $n \in \mathbb{Z}$ tels que  $D(a,n) \cap \varpi_F\OF = \emptyset$, tout $0\leqslant \underline{n}_{S\setminus J}\leqslant \underline{d}_{S\setminus J}$ et tout $\underline{m}_J \in \mathbb{Z}_{\geq 0}^{|J|}$;
\item[(ii)] \eqref{puno} pour $a=0$, pour tout entier $n \leq 0$, tout $0\leqslant \underline{n}_{S\setminus J}\leqslant \underline{d}_{S\setminus J}$ et tout $\underline{m}_J \in \mathbb{Z}_{\geq 0}^{|J|}$  tels que $r- (|\underline{n}_{S\setminus J}|+|\underline{m}_J|)>0$;
\item[(iii)] \eqref{pdue}  pour tout $a \in F$, tout $n\geq 0$, tout $\underline{0}\leqslant \underline{n}_{S\setminus J} \leqslant \underline{d}_{S\setminus J}$ et tout $\underline{m}_J \in \mathbb{Z}_{\geq 0}^{|J|}$ tels que $r-(|\underline{n}_{S\setminus J}|+|\underline{m}_J|) > 0$.
\end{itemize}

La proposition suivante montre que \eqref{trecrit} implique (i). 

\begin{prop}\label{hhhh} 
La condition \eqref{trecrit} implique la condition \eqref{puno} pour tout disque $D(a,n)$ avec $a\in F$ et $n \in \mathbb{Z}$ tel que $D(a,n) \cap \varpi_F \OF = \emptyset$, tout $0\leqslant \underline{n}_{S\setminus J}\leqslant \underline{d}_{S\setminus J}$ et tout $\underline{m}_J \in \mathbb{Z}_{\geq 0}^{|J|}$. 
\begin{proof}

Un calcul analogue à celui de la Proposition \ref{messi1} et dont on laisse les détails au lecteur, montre que la condition \eqref{trecrit} est équivalente à
\begin{align} 
&\Big|\int_{D(\frac{1}{a},n-\frac{2 val_F(a)}{f})}  z^{ \underline{d}_{S\setminus J}} \Big(\frac{1}{z}-a\Big)^{\underline{n}_{S\setminus J}} \Big(\frac{1}{z}-a\Big)^{\underline{m}_J} \mu(z) \Big| \leq C_{\mu} |a|^{2r - |\underline{d}_{S\setminus J}|} q^{n (r-|\underline{n}_{S\setminus J}|-|\underline{m}_J|)} \label{quattrocrit} 
\end{align}
pour tout $a \in \OF-\{0 \}$, tout $\underline{0} \leqslant \underline{n}_{S\setminus J}\leqslant \underline{d}_{S\setminus J}$, tout $\underline{m}_J \in \mathbb{Z}_{\geq 0}^{|J|}$ et tout entier $n> \frac{val_F(a)}{f}$.

Soit $a \in F$ et $n \in \mathbb{Z}$ tel que $D(a,n) \cap \varpi_F \OF = \emptyset$. Pour tout  $0\leqslant \underline{n}_{S\setminus J}\leqslant \underline{d}_{S\setminus J}$ on a les identités suivantes:  
\begin{align*}
\mathbf{1}_D(z) \Big(z - \frac{1}{a}  \Big)^{\underline{n}_{S\setminus J}} &= \mathbf{1}_D(z)(-1)^{\underline{n}_{S\setminus J}} a^{-\underline{n}_{S\setminus J}} z^{\underline{n}_{S\setminus J}} \Big(\frac{1}{z} - a \Big)^{\underline{n}_{S\setminus J}} \\
&= \mathbf{1}_D(z)(-1)^{\underline{n}_{S\setminus J}} a^{-\underline{n}_{S\setminus J}} \Big(\frac{1}{z} - a + a \Big)^{\underline{d}_{S\setminus J} -\underline{n}_{S\setminus J}} z^{\underline{d}_{S\setminus J}} \Big(\frac{1}{z} - a \Big)^{\underline{n}_{S\setminus J}} \\
&= \mathbf{1}_D(z) \sum_{\underline{0}\leqslant \underline{k}_{S\setminus J}\leqslant \underline{d}_{S\setminus J} -\underline{n}_{S\setminus J}}   \lambda_{\underline{k}_{S\setminus J}}a^{\underline{k}_{S\setminus J}-\underline{n}_{S\setminus J}}  z^{\underline{d}_{S\setminus J}} \Big(\frac{1}{z}-a  \Big)^{\underline{d}_{S\setminus J} - \underline{k}_{S\setminus J}}   
\end{align*}
où les $\lambda_{\underline{k}_{S\setminus J}}$ sont des entiers. Par un calcul similaire au précédent on obtient:
\[
\mathbf{1}_D(z) \Big(z - \frac{1}{a}  \Big)^{\underline{m}_{J}} = \mathbf{1}_D(z) \sum_{ \underline{r}_{J}\geqslant \underline{0}}   \mu_{\underline{r}_{J}} a^{-2\underline{m}_{ J}-\underline{r}_{J}} \Big(\frac{1}{z}-a  \Big)^{\underline{r}_{J} + \underline{m}_{J}} 
\]
où les $\mu_{\underline{r}_{J}}$ sont des entiers.

Les deux identités ci-dessus et la condition \eqref{quattrocrit} impliquent:
\begin{align*}
& \,  \Big| \mu \Big( \mathbf{1}_{D}(z) \Big( z-\frac{1}{a} \Big)^{\underline{n}_{S\setminus J}} ( z-\frac{1}{a} \Big)^{\underline{m}_{J}}   \Big)\Big| \\
\leq & \, \Big| a^{-\underline{n}_{S\setminus J}}a^{-2\underline{m}_{J}}   \sup_{\substack{
\underline{r}_{J}\geqslant \underline{0}\\
\underline{0}\leqslant \underline{k}_{S\setminus J}\leqslant \underline{d}_{S\setminus J}-\underline{n}_{S\setminus J}
}}   a^{\underline{k}_{S\setminus J}} a^{-\underline{r}_{J}} \mu\Big( \mathbf{1}_{D}(z) z^{\underline{d}_{S\setminus J}} \Big(\frac{1}{z}-a\Big)^{\underline{d}_{S\setminus J}-\underline{k}_{S\setminus J}} \Big(\frac{1}{z}-a\Big)^{\underline{r}_{J}+\underline{m}_{J}} \Big)   \Big| \\
\leq & \, C_{\mu} |a|^{-|\underline{n}_{S\setminus J}|-2|\underline{m}_{J}|}   \sup_{\substack{\underline{r}_{J}\geqslant \underline{0} \\
\underline{0}\leqslant \underline{k}_{S\setminus J}\leqslant \underline{d}_{S\setminus J}-\underline{n}_{S\setminus J}
}}   |a|^{|\underline{k}_{S\setminus J}|-|\underline{r}_{J}|}   |a|^{2r-|\underline{d}_{S\setminus J}|}    q^{n(r - |\underline{d}_{S\setminus J}|+|\underline{k}_{S\setminus J}|-|\underline{r}_{J}| -|\underline{m}_{J}|)} \\
= &\, C_{\mu} |a|^{2r-2|\underline{n}_{S\setminus J}|-2|\underline{m}_{J}|} q^{n(r - |\underline{n}_{S\setminus J}| -|\underline{m}_{J}|)}  \\
= & \, C_{\mu} q^{(n- 2\frac{val_F(a)}{f})(r- |\underline{n}_{S\setminus J}|-|\underline{m}_J|)}.
\end{align*}
Comme $D\big(\frac{1}{a}, n-\frac{val_F(a)}{f}\big)$ pour $a \in \OF-\{0 \}$ et $n > \frac{val_F(a)}{f}$ parcourt tous les disques $D(b,m)$ avec $b\in F$ et $m\in \mathbb{Z}_{\geq 0}$ dans $F$ tels que $D(b,m) \cap \varpi_F \OF = \emptyset$, on peut conclure. 
\end{proof}
\end{prop}

En utilisant les conditions \eqref{funo} et \eqref{fdue} on voit que montrer $(ii)$ et $(iii)$ est équivalente à montrer (quitte à modifier la constante $C_{\mu}$)
\begin{align}\label{altracond1}
&\Big|\int_{F\setminus D(0,n)} z^{\underline{n}_{S\setminus J}} z^{\underline{m}_J} \mu(z)\Big| \leq C_{\mu} q^{n (r-|\underline{n}_{S\setminus J}|-|\underline{m}_J|)}
\end{align}
pour tout entier $n \leq 0$, tout $0\leqslant \underline{n}_{S\setminus J}\leqslant \underline{d}_{S\setminus J}$ et tout $\underline{m}_J \in \mathbb{Z}_{\geq 0}^{|J|}$  tels que $r- (|\underline{n}_{S\setminus J}|+|\underline{m}_J|)>0$ et
\begin{align}\label{altracond2}
&\Big|\int_{D(a,n+1)} \chi_2\chi_1^{-1}(z-a) (z-a)^{\underline{d}_{S\setminus J}- \underline{n}_{S\setminus J}} (z-a)^{-\underline{m}_J} \mu(z)\Big| \leq C_{\mu} q^{n(|\underline{n}_{S\setminus J}|+|\underline{m}_{J}|-r)}
\end{align}
pour tout $a \in F$, tout $n\geq 0$, tout $\underline{0}\leqslant \underline{n}_{S\setminus J} \leqslant \underline{d}_{S\setminus J}$ et tout $\underline{m}_J \in \mathbb{Z}_{\geq 0}^{|J|}$ tels que $r-(|\underline{n}_{S\setminus J}|+|\underline{m}_J|) > 0$.

Rappelons que si $f \in B(\chi,J,\underline{d}_{S\setminus J})$ alors
\begin{align}\label{normabb}
\|f \|_{B} = \sup \big(\|f_1\|_{C^r}, \|f_2\|_{C^r}\big)
\end{align}
où $(f_1,f_2)$ désigne l'élément de $C^{r}\big(\OF,J,\underline{d}_{S\setminus J}\big)^2$ qui correspond à $f$ via l'isomorphisme \eqref{isocr}. 

Les conditions \eqref{altracond1} et \eqref{altracond2} sont une conséquence immédiate du lemme suivant. 

\begin{lemm}\label{lemmafinale1}
\begin{itemize} 
\item[$\bullet$] Il existe une constante $C \in \mathbb{R}_{\geq 0}$ telle que pour tout entier $n \leq 0$, tout $\underline{0}\leqslant \underline{n}_{S\setminus J}\leqslant \underline{d}_{S\setminus J}$ et tout $\underline{m}_J \in \mathbb{Z}_{\geq 0}^{|J|}$  tels que $r- (|\underline{n}_{S\setminus J}|+|\underline{m}_J|)>0$ on a:
\begin{align*}
\|\mathbf{1}_{F\setminus D(0,n+1)}(z) z^{\underline{n}_{S\setminus J}} z^{\underline{m}_J}  \|_{B} \leq C q^{n (r-|\underline{n}_{S\setminus J}|-|\underline{m}_J|)}.
\end{align*}
\item[$\bullet$] Il existe une constante $C \in \mathbb{R}_{\geq 0}$ telle que pour tout $a \in F$, tout $n\geq 1$, tout $\underline{0}\leqslant \underline{n}_{S\setminus J} \leqslant \underline{d}_{S\setminus J}$ et tout $\underline{m}_J \in \mathbb{Z}_{\geq 0}^{|J|}$ tels que $r-(|\underline{n}_{S\setminus J}|+|\underline{m}_J|) > 0$ on a:
\begin{align*}
\| \mathbf{1}_{D(a,n)}(z) \chi_2\chi_1^{-1}(z-a) (z-a)^{\underline{d}_{S\setminus J}- \underline{n}_{S\setminus J}} (z-a)^{-\underline{m}_J} \|_{B} \leq C q^{n(|\underline{n}_{S\setminus J}|+|\underline{m}_{J}|-r)}.
\end{align*}
\end{itemize}
\begin{proof}
Notons $f_{\underline{n}_{S\setminus J},\underline{m}_J}$, pour tout $\underline{0}\leqslant \underline{n}_{S\setminus J}\leqslant \underline{d}_{S\setminus J}$ et tout $\underline{m}_J \in \mathbb{Z}_{\geq 0}^{|J|}$ tels que $r-(|\underline{n}_{S\setminus J}|+|\underline{m}_J|) > 0$,  la fonction de $\OF$ dans $E$ définie par:  
\[
z \mapsto \chi_2 \chi_1^{-1}(z) z^{\underline{d}_{S\setminus J}-\underline{n}_{S\setminus J}}z^{-\underline{m}_{J}}. 
\] 
D'après le Lemme \ref{funzionicr} c'est une fonction de classe $C^r$. Posons:
\begin{align}
C = \sup \Big\{ \|f_{\underline{n}_{S\setminus J},\underline{m}_J} \|_{C^r}: \underline{0}\leqslant \underline{n}_{S\setminus J}\leqslant \underline{d}_{S\setminus J},\, \underline{m}_J \in \mathbb{Z}_{\geq 0}^{|J|} \ \mbox{et}\ r- (|\underline{n}_{S\setminus J}|+|\underline{m}_J|)>0  \Big\}.
\end{align}

\begin{itemize}
\item[$\bullet$] Par \eqref{normabb} on a:
\begin{align*}
 \|\mathbf{1}_{F\setminus D(0,n+1)}(z) z^{\underline{n}_{S\setminus J}} z^{\underline{m}_J}  \|_{B} = \|\mathbf{1}_{D(0,-n)}(z) f_{\underline{n}_{S\setminus J},\underline{m}_J}(z) \|_{C^r}. \\
\end{align*}
On peut récrire $\|\mathbf{1}_{D(0,-n)}(z) f_{\underline{n}_{S\setminus J},\underline{m}_J}(z) \|_{C^r}$ sous la forme:
\[
 \ \Big|\chi_2 \chi_1^{-1}(\uni^{-n}) (\uni^{-n})^{\underline{d}_{S\setminus J} - \underline{n}_{S\setminus J}} (\uni^{-n})^{-\underline{m}_{J}}  \Big|   \Big\|\mathbf{1}_{D(0,-n)}(z) f_{\underline{n}_{S\setminus J},\underline{m}_J}\Big( \frac{z}{\uni^{-n}} \Big)  \Big\|_{C^r}. 
\]
Or, d'après \eqref{carattereintegro} on a:
\[
\Big|\chi_2 \chi_1^{-1}(\uni^{-n}) (\uni^{-n})^{\underline{d}_{S\setminus J} - \underline{n}_{S\setminus J}} (\uni^{-n})^{-\underline{m}_{J}}  \Big| = q^{n(2r-|\underline{n}_{S\setminus J}|-|\underline{m}_J|)}
\] 
et d'après le Lemme \ref{corollarioutile} on a:
\[
\Big\|\mathbf{1}_{D(0,-n)}(z) f_{\underline{n}_{S\setminus J},\underline{m}_J}\Big( \frac{z}{\uni^{-n}} \Big)  \Big\|_{C^r} \leq C q^{-nr},
\]
d'où
\[
\|\mathbf{1}_{F\setminus D(0,n+1)}(z) z^{\underline{n}_{S\setminus J}} z^{\underline{m}_J}  \|_{B} \leq C q^{n(r-|\underline{n}_{S\setminus J}|-|\underline{m}_J|)}.
\]

\item[$\bullet$] On distingue deux cas.

(i) Supposons $a \in \varpi_F \OF$. Par \eqref{normabb} on a:
\begin{align*}
\| \mathbf{1}_{D(a,n)}(z) \chi_2\chi_1^{-1}(z-a) (z-a)^{\underline{d}_{S\setminus J}- \underline{n}_{S\setminus J}} &(z-a)^{-\underline{m}_J} \|_{B} \\ 
 &= \| \mathbf{1}_{D(\frac{a}{\varpi_F},n-1)}(z) f_{\underline{n}_{S\setminus J},\underline{m}_J}(\varpi_F z-a) \|_{C^r},
\end{align*}
et comme la norme $C^r$ est invariante par translation on déduit l'égalité suivante:
\[
\| \mathbf{1}_{D(\frac{a}{\varpi_F},n-1)}(z) f_{\underline{n}_{S\setminus J},\underline{m}_J}(\varpi_F z-a) \|_{C^r} = \| \mathbf{1}_{D(0,n-1)}(z) f_{\underline{n}_{S\setminus J},\underline{m}_J}(\varpi_F z) \|_{C^r}.
\]
On peut récrire $\| \mathbf{1}_{D(0,n-1)}(z) f_{\underline{n}_{S\setminus J},\underline{m}_J}(\varpi_F z) \|_{C^r}$ sous la forme:
\[
\Big|\chi_2\chi_1^{-1}(\varpi_F^n) (\varpi_F^n)^{(\underline{d}_{S\setminus J}-\underline{n}_{S\setminus J})} (\varpi_F^n)^{-\underline{m}_J}  \Big| \Big\|\mathbf{1}_{D(0,n-1)}(z) f_{\underline{n}_{S\setminus J},\underline{m}_J} \Big(\frac{z}{\varpi_F^{n-1}} \Big)  \Big\|_{C^r}.
\]
D'après \eqref{carattereintegro} on a:
\[
\Big|\chi_2\chi_1^{-1}(\varpi_F^n) (\varpi_F^n)^{(\underline{d}_{S\setminus J}-\underline{n}_{S\setminus J})} (\varpi_F^n)^{-\underline{m}_J}  \Big| = q^{n(-2r + |\underline{n}_{S\setminus J}|+ |\underline{m}_{J}|)}
\]
et d'après le Lemme \ref{corollarioutile} on a:
\[
 \Big\|\mathbf{1}_{D(0,n-1)}(z) f_{\underline{n}_{S\setminus J},\underline{m}_J} \Big(\frac{z}{\varpi_F^{n-1}} \Big)  \Big\|_{C^r} \leq C q^{(n-1)r}
\]
d'où
\[
\| \mathbf{1}_{D(a,n)}(z) \chi_2\chi_1^{-1}(z-a) (z-a)^{\underline{d}_{S\setminus J}- \underline{n}_{S\setminus J}} (z-a)^{-\underline{m}_J} \|_{B} \leq C q^{-r} q^{n(-r + |\underline{n}_{S\setminus J}|+ |\underline{m}_{J}|)}.
\]

(ii) Supposons $a \notin \varpi_F\OF$. Par \eqref{normabb} on a:
\begin{align*}
\quad \quad \| &\mathbf{1}_{D(a,n)}(z) \chi_2\chi_1^{-1}(z-a) (z-a)^{\underline{d}_{S\setminus J}- \underline{n}_{S\setminus J}} (z-a)^{-\underline{m}_J} \|_{B} \\
 &\quad = \Big|\chi_2\chi_1^{-1}(a)a^{\underline{d}_{S\setminus J}-\underline{n}_{S\setminus J}} a^{-\underline{m}_{J}} \Big|  \Big\|\mathbf{1}_{D\big(\frac{1}{a},n-\frac{2 val_F(a)}{f}\big)}(z) z^{\underline{n}_{S\setminus J}} z^{\underline{m}_J} f_{\underline{n}_{S\setminus J},\underline{m}_J}\Big(z-\frac{1}{a}  \Big)  \Big\|_{C^r}.
\end{align*}
En écrivant $z^{\underline{n}_{S\setminus J}} = (z-\frac{1}{a}+\frac{1}{a})^{\underline{n}_{S\setminus J}}$ (resp. $z^{\underline{m}_{J}} = (z-\frac{1}{a}+\frac{1}{a})^{\underline{m}_{J}}$) et en développant, et en utilisant l'invariance par translation de la norme $C^r$ on obtient:
\begin{align*}
\quad\quad \Big\|\mathbf{1}_{D\big(\frac{1}{a},n-\frac{2 val_F(a)}{f}\big)}(z) & z^{\underline{n}_{S\setminus J}} z^{\underline{m}_J} f_{\underline{n}_{S\setminus J},\underline{m}_J}\Big(z-\frac{1}{a}  \Big)  \Big\|_{C^r} \\ 
&\quad \leq \sup_{\substack{
\underline{0} \leqslant \underline{\alpha}_J \leqslant \underline{m}_J  \\
\underline{0} \leqslant \underline{\beta}_{S\setminus J} \leqslant \underline{n}_{S\setminus J} }}  |a|^{-|\underline{\alpha}_J|-|\underline{\beta}_{S\setminus J}|} \Big\| \mathbf{1}_{D\big(0,n-\frac{2 val_F(a)}{f}\big)}(z) f_{\underline{\beta}_{S\setminus J},\underline{\alpha}_J}(z) \Big\|_{C^r}.
\end{align*} 
Par le Lemme \ref{corollarioutile} on a:
\begin{align*}
\quad\quad \Big\| \mathbf{1}_{D\big(0,n-\frac{2 val_F(a)}{f}\big)}(z) f_{\underline{\beta}_{S\setminus J},\underline{\alpha}_J}&(z) \Big\|_{C^r} \\ 
&\quad \leq C \Big|\chi_2\chi_1^{-1}\Big(\frac{\varpi_F^n}{a^2} \Big)  \Big(\frac{\varpi_F^n}{a^2} \Big)^{\underline{d}_{S\setminus J}-\underline{\beta}_{S\setminus J}} \Big(\frac{\varpi_F^n}{a^2} \Big)^{-\underline{\alpha}_J}  \Big| \Big| \frac{\varpi_F^n}{a^2} \Big|^{-r}, 
\end{align*}
et comme le $\sup$ est atteint pour $\underline{\alpha}_J = \underline{m}_J$, $\underline{\beta}_{S\setminus J} = \underline{n}_{S\setminus J}$ on obtient, en utilisant \eqref{carattereintegro}: 
\[
\| \mathbf{1}_{D(a,n)}(z) \chi_2\chi_1^{-1}(z-a) (z-a)^{\underline{d}_{S\setminus J}- \underline{n}_{S\setminus J}} (z-a)^{-\underline{m}_J} \|_{B} \leq C q^{r  - |\underline{n}_{S\setminus J}| - |\underline{m}_J|},
\]
d'où le résultat.
\end{itemize}
\end{proof}
\end{lemm}

Le Lemme \ref{lemmafinale1} termine la preuve de $(B)\Rightarrow (A)$, et donc la preuve du Théorème \ref{teoremaprincipale}. Ainsi on a obtenu que l'espace de Banach dual du complété cherché est isomorphe dans $I(\chi,J,\underline{d}_{S\setminus J})^{\vee}$ au sous-espace de Banach de $B(\chi,J,\underline{d}_{S\setminus J})^{\vee}$ formé des $\mu$ qui annulent $L(\chi,J,\underline{d}_{S\setminus J})$, c'est-à-dire à $\Pi(\chi,J,\underline{d}_{S\setminus J})^{\vee}$. En particulier $\Pi(\chi,J,\underline{d}_{S\setminus J})^{\vee}$ est un $G$-Banach unitaire. \\

Rappelons que dans \cite{sch3} est introduite la catégorie $\mathrm{Mod}_{comp}^{fl}(\Oe)$ des $\Oe$-modules sans torsion linéairement topologiques séparés compacts, les morphismes étant les applications $\Oe$-linéaires continues. Soit $M \in \mathrm{Mod}_{comp}^{fl}(\Oe)$ et définissons le $E$-espace de Banach:
\[
M^d \ugu \mathrm{Hom}_{\Oe}^{cont}(M,E) \ \mbox{muni de la norme} \ \|l\| \ugu \sup_{x \in M} |l(x)|.
\]  
Notons $\mathrm{Mod}_{comp}^{fl}(\Oe)_{\mathbb{Q}}$ la catégorie ayant les mêmes objets que la catégorie  $\mathrm{Mod}_{comp}^{fl}(\Oe)$ mais pour morphismes:
\[
\mathrm{Hom}_{\mathrm{Mod}_{comp}^{fl}(\Oe)_{\mathbb{Q}}}(A,B) \ugu \mathrm{Hom}_{\mathrm{Mod}_{comp}^{fl}(\Oe)}(A,B) \otimes E.
\] Dans \cite[Théorème 1.2]{sch3}, il est montré que le foncteur $M \mapsto M^d$ induit une anti-équivalence de catégories entre $\mathrm{Mod}_{comp}^{fl}(\Oe)_{\mathbb{Q}}$ et la catégorie des $E$-espaces de Banach.

\begin{coro}\label{corfin}
Il existe un isomorphisme $G$-équivariant d'espaces de Banach $p$-adiques: 
\[
I(\chi,J,\underline{d}_{S\setminus J})^{\bigwedge} \stackrel{\sim}{\longrightarrow} \Pi(\chi,J,\underline{d}_{S\setminus J}).
\]
\begin{proof}
L'argument est analogue à celui donné dans \cite[Théorème 4.3.1]{bb}. D'après \cite[Lemme 9.9]{sch} on a une injection fermée $G$-équivariante:
\[
\Pi(\chi,J,\underline{d}_{S\setminus J}) \into \Big(\Pi(\chi,J,\underline{d}_{S\setminus J})^{\vee}\Big)^{\vee}.
\]
Cela implique que $\Pi(\chi,J,\underline{d}_{S\setminus J})$ est aussi un $G$-Banach unitaire. Alors, par la propriété universelle du complété unitaire universel, l'application $I(\chi,J,\underline{d}_{S\setminus J}) \to \Pi(\chi,J,\underline{d}_{S\setminus J})$ induit un morphisme $G$-équivariant continu de $I(\chi,J,\underline{d}_{S\setminus J})^{\bigwedge}$ vers $\Pi(\chi,J,\underline{d}_{S\setminus J})$. Cela induit un morphisme continu sur les duaux munis de leur topologie faible qui sont des éléments de $\mathrm{Mod}_{comp}^{fl}(\Oe)_{\mathbb{Q}}$. Or, par le Théorème \ref{teoremaprincipale} ce morphisme est bijectif et continu. Donc, d'après  \cite[Lemme 4.2.2]{breuila} c'est aussi un isomorphisme pour les topologies faibles. Par dualité (\cite[Théorème 1.2]{sch3}) on obtient alors l'isomorphisme topologique $\GL(F)$-équivariant de l'énoncé. 
\end{proof}
\end{coro}

\begin{rema}{\rm Le Corollaire \ref{corfin} généralise le \cite[Théorème 4.3.1]{bb} pour $F=\Q$. Mentionnons que ce résultat joue un rôle important dans la preuve par Berger et Breuil de la non nullité de l'espace $I(\chi,J,\underline{d}_{S\setminus J})^{\bigwedge}$.}
\end{rema}

\subsection{Exemple}

Introduisons quelque notations supplémentaires et rappelons la construction des représentations considérées dans \cite{bre}. Si $\lambda \in E^{\times}$ on désigne par $\mathrm{unr}_F(\lambda)\colon F^{\times} \to E^{\times}$ le caractère non ramifié défini par $x \mapsto \lambda^{val_F(x)}$. Soient $\alpha,\tilde{\alpha} \in E^{\times}$ et $\underline{k} \in \mathbb{Z}_{>1}^{|S|}$. Fixons $J_1,J_2$ deux sous-ensembles de $S$ tels que $J_1\subseteq J_2 \subseteq S$. Considérons les deux caractères algébriques suivants:
\[
\chi_1 = \mathrm{unr}_F(\alpha^{-1}) \prod_{\sigma \in J_1} \sigma^{k_{\sigma}-1}, \quad \chi_2 = \mathrm{unr}_F(p\tilde{\alpha}^{-1}) \prod_{\sigma \in J_1} \sigma^{-1} \prod_{\sigma \in J_2\setminus J_1} \sigma^{k_{\sigma}-2} 
\]
et posons:
\[
\pi(J_1,J_2) = \Big(  \bigotimes_{\sigma \in S \setminus J_2} (\Sym^{k_{\sigma}-2} E^2)^{\sigma}  \Big) \otimes_E \Big(\Ind_{P}^{G} \chi_1 \otimes \chi_2 \Big)^{J_2-{an}}. 
\]
D'après la Proposition \ref{nullita} on connait deux conditions nécessaires pour que le complété unitaire universaire de la représentation $\Q$-analytique $\pi(J_1,J_2)$ soit non nul. Un calcul immédiat montre qu'elles sont équivalentes à
\begin{align}
-(\val_{F}(\alpha)+ \val_{F}(\tilde{\alpha})) + \sum_{\sigma\in S} (k_{\sigma}-1) &= 0 \label{unoesempio} \\
-\val_{F}(\tilde{\alpha}) + \sum_{\sigma \in S\setminus J_1}(k_{\sigma}-1) &\geq 0. \label{dueesempio}
\end{align} 
Supposons que \eqref{unoesempio} et \eqref{dueesempio} soient satisfaits. En particulier on en déduit l'inégalité suivante  
\[
-val_F(\alpha)+\sum_{\sigma \in J_1}(k_{\sigma}-1)\leq 0.
\]
Notons $r = val_F(\alpha)-\sum_{\sigma \in J_1}(k_{\sigma}-1)$ et 
\[
J_3 = J_2 \coprod \{\sigma \in S\backslash J_2,\, k_{\sigma}-1 > r  \}.
\]
D'après la Proposition \ref{dsigma} on sait que l'application fermée et $G$-équivariante
\[
\pi(J_1,J_2) \into \pi(J_1,J_3) \ugu \Big(  \bigotimes_{\sigma \in S \setminus J_3} (\Sym^{k_{\sigma}-2} E^2)^{\sigma}  \Big) \otimes_E \Big(\Ind_{P}^{G} \chi_1 \otimes \chi_2 \prod_{\sigma \in J_3\setminus J_2} \sigma^{k_{\sigma}-2} \Big)^{J_3-{an}}
\]
induit un isomorphisme $G$-équivariante de $\pi(J_1,J_2)^{\bigwedge}$ dans $\pi(J_1,J_3)^{\bigwedge}$. Posons:
\[
\chi_1' = \chi_1, \quad \chi_2' =  \chi_2 \prod_{\sigma \in J_3\setminus J_2} \sigma^{k_{\sigma}-2}. 
\]
Considérons:
\[
B(\chi, J_3, (k_{\sigma}-2)_{\sigma\notin J_3}) = C^r(\OF, J_3, (k_{\sigma}-2)_{\sigma\notin J_3}) \oplus C^r(\OF, J_3, (k_{\sigma}-2)_{\sigma\notin J_3}).
\]
C'est un espace de Banach sur $E$ muni d'une action continue de $G$ (voir la preuve du Lemme \ref{autocont}). D'après le Lemme \ref{funzionicr} la fonction $h_{(n_{\sigma})_{\sigma \notin  J_3},(m_{\sigma})_{\sigma\in J_3}}$ définie par:
\[
z \mapsto \chi_2' {\chi_1'}^{-1}(z) \prod_{\sigma\notin J_3} \sigma(z)^{k_{\sigma}-2-n_{\sigma}} \prod_{\sigma \in J_3} \sigma(z)^{-m_{\sigma}}
\]
se prolonge sur $\OF$ en une fonction de classe $C^r$. On désigne par $L(\chi, J_3,(k_{\sigma}-2)_{\sigma\notin J_3})$ le sous-espace de $B(\chi, J_3, (k_{\sigma}-2)_{\sigma\notin J_3})$ engendré par les couples de fonctions:
\begin{align*}
&\Big(z\mapsto \prod_{\sigma\notin J_3}\sigma(\uni z)^{n_{\sigma}} \prod_{\sigma\in J_3}\sigma(\uni z)^{m_{\sigma}}, z\mapsto h_{(n_{\sigma})_{\sigma \notin  J_3},(m_{\sigma})_{\sigma\in J_3}}(z) \Big) \\
&\Big(z\mapsto h_{(n_{\sigma})_{\sigma\notin J_3},(m_{\sigma})_{\sigma \in J_3}}(\uni z -a), z\mapsto h_{(n_{\sigma})_{\sigma\notin J_3},(m_{\sigma})_{\sigma \in J_3}}(1-az) \prod_{\sigma \notin  J_3} \sigma(z)^{n_{\sigma}}  \prod_{\sigma \in J_3} \sigma(z)^{m_{\sigma}}  \Big)
\end{align*}
pour tout $a\in F$, tout $(m_{\sigma})_{\sigma \in J_3} \in \mathbb{Z}_{\geq 0}^{|J_3|}$ et tout $(n_{\sigma})_{\sigma \notin J_3} \leq (k_{\sigma}-2)_{\sigma \notin J_3}$ tels que  $r - \sum_{\sigma \notin J_3} n_{\sigma} - \sum_{J_3} m_{\sigma} > 0$. Alors par le Corollaire \ref{corfin} on a:  
\[
\pi(J_1,J_2)^{\bigwedge} \stackrel{\sim}{\longrightarrow} B(\chi, J_3, \underline{k}_{S\setminus J_3}-2)/L(\chi, J_3, \underline{k}_{S\setminus J_3}-2).  
\]




\begin{thebibliography}{1}

\bibitem{ami}
Y. Amice,  {\it  Duals}, Proc. of a conf. on $p$-adic analysis (Nijmegen 1978), Nijmegen, Math. Institut Katholische Univ., 1978, 1-15.

\bibitem{amivel}
Y. Amice \& J. Vélu,  {\it  Distributions $p$-adiques associées aux séries de Hecke}, Astérisque 24-25, 1975, 119-131.


\bibitem{bb}
L. Berger et C. Breuil,  {\it  Sur quelques représentations potentiellement cristallines de $\mathrm{GL}_2(\mathbb{Q}_p)$}, Astérisque 330, 2010, 155-211.

\bibitem{breuilab}
C. Breuil,  {\it  Sur quelques représentations modulaires et $p$-adique de $\GL(\Q)$ II}, J. Inst. Math. Jussieu  2, 2003, 1-36.

\bibitem{breuila}
C. Breuil,  {\it  Invariant $L$ et série spéciale $p$-adique}, Ann. Scient. de l'E.N.S  37, 2004, 559-610.

\bibitem{breuilb}
C. Breuil,  {\it  Série spéciale $p$-adique et cohomologie étale complétée}, Astérisque 331, 2010, 65-115.


\bibitem{bre}
C. Breuil,  {\it  Remarks on some locally $\Q$-analytic representations of $\GL(F)$ in the crystalline case}, Non-abelian Fundamental Groups and Iwasawa Theory, Cambridge University Press, 2012.


\bibitem{be}
C. Breuil et M. Emerton,  {\it  Représentations $p$-adiques ordinaires de $\GL(\Q)$ et compatibilité local-global}, Astérisque 331, 2010, 255-315.

\bibitem{bs} 
C. Breuil et P. Schneider,  {\it First steps towards $p$-adic Langlands functoriality}, J. Reine Angew. Math. 610, 2007, 149-180.

\bibitem{bv}
N. Bourbaki,  {\it  Vari\'{e}tés différentielles et analytiques. Fascicule de résultats}, Paris: Hermann 1967.


\bibitem{colmez}
P. Colmez,  {\it Fonctions d'une variable $p$-adique}, Astérisque 330, 2010, 13-59.

\bibitem{colmez3}
P. Colmez,  {\it La série principale unitaire de $\GL(\Q)$}, Astérisque  330, 2010, 213-262.


\bibitem{colmez2}
P. Colmez,  {\it Représentations de $\GL(\Q)$ et $(\varphi,\Gamma)$-modules}, Astérisque 330, 2010, 281-509.

\bibitem{sha}
E. De Shalit,  {\it Mahler bases, lubin tate groups and elementary $p$-adic analysis}, Prépublication.

\bibitem{marco}
M. De Ieso,  {\it Espaces de fonctions de classe $C^r$ sur $\OF$}, Prépublication.

\bibitem{marco2}
M. De Ieso,  {\it Un critère de non nullité}, En préparation.



\bibitem{eme}
M. Emerton,  {\it  $p$-adic $L$-functions and unitary completions of representations of $p$-adic reductive groups}, Duke Math. J. 130, 2005, 353-392.


\bibitem{font}
J.-M. Fontaine,  {\it  Représentations $p$-adiques de corps locaux I}, The Grothendieck Festschrift, Vol II, Progr. Math., Birkhauser, vol 87, 1990, 249-309.




\bibitem{fr}
H. Frommer,  {\it  The locally analyitic principal series of split reductive groups}, Prépublication 2003.




 
\bibitem{ks}
D. Kazhdan et E. Shalit  {\it Kirillov models and integrals structures in $p$-adic smooth representations of $\GL(F)$ }, preprint. 

\bibitem{enno1}
E. Nagel,  {\it  Fractional non-archimedean calculus in one variable}, Prépublication.


\bibitem{enno2}
E. Nagel,  {\it  Fractional non-archimedean calculus in many variables}, Prépublication.





\bibitem{pas}
V. Pa\v{s}k\={u}nas,  {\it Admissible unitary completions of locally $\Q$-rational representations of $\GL(F)$}, Representation theory 14, 2010, 324-354.  

\bibitem{pask}
V. Pa\v{s}k\={u}nas,  {\it The image of Colmez's Montreal Functor}, prépublication, 2011. 
 

\bibitem{sch}
P. Schneider, {\it Nonarchimedean Functional Analysis}, Springer Monographs in Mathematics, Springer Verlag, 2002.

\bibitem{sch1}
P. Schneider et J. Teitelbaum, {\it Locally analytic distributions and $p$-adic representation theory, with an application to $\GL$}, J. Amer. Math. Soc. 15, 2002, 51-125.

\bibitem{sch3}
P. Schneider et J. Teitelbaum, {\it Banach space representations and Iwasawa theory}, Israel J. Math. 127, 2002, 359-380.




\bibitem{scr}
B. Schraen, {\it Représentations $p$-adiques de $\GL(F)$ et catégories dérivées}, à paraître dans Israel Journal of Math. 176.

\bibitem{vander}
M. Van der Put, {\it Algèbres de fonctions continues $p$-adiques}, Proc. Kon. Ned. Akad. v. Wetensch. A 71, 1968, 556-661.


\bibitem{vig}
M.-F. Vigneras, {\it A criterion for integral structures and coefficients systems on the three of $PGL(2,F)$ }, Pure and Applied Mathematics Quaterly 4, 2008, 1291-1316.

\bibitem{vis}
M. Vishik, {\it Non-archimedian measures connected with Dirichlet series}, Math. USSR Sbornik   28, 1976, 216-228.


\end{thebibliography}
\end{document}